\documentclass[11pt,oneside,reqno]{amsart}
\usepackage{amssymb,amsmath,amsthm,amsopn,amstext,amsfonts,amsbsy}
\usepackage[english]{babel}
\font\ninerm=cmr9
\font\ninebf=cmbx9

\title[Transition asymptotics for random walks]{
Transition asymptotics for reaction-diffusion in random media}

\author{G\'erard Ben Arous$^1$, Stanislav Molchanov$^1$ and Alejandro F. Ram\'\i rez$^{1,2}$}

\thanks{ AMS 2000 {\it subject classifications}. Primary  82B41, 82B44 ;
secondary 60J45, 60J65, 82C22.}

\thanks{$^1$Partially supported by Fondo Nacional de Desarrollo Cient\'\i fico
y Tecnol\'ogico grant 7020686}

\thanks{$^2$Partially supported by Fondo Nacional de Desarrollo Cient\'\i fico
y Tecnol\'ogico grant 1020686}

\thanks{{\it Key words and phrases.} Parabolic Anderson model,
random walk, branching processes, Feynman-Kac formula, principal eigenvalue.}

\address[G\'erard Ben Arous]{Ecole Polytechnique F\'ed\'erale de Lausanne, CH-1015
Switzerland\\
\and Courant Institute of Mathematical Sciences\\ New York University\\
 251 Mercer Street, NY 10012, USA}

\address[Stanislav Molchanov]{Department of Mathematics\\
University of North Carolina-Charlotte\\
376 Fretwell Bldg. 9201 University City Blvd.\\
Charlotte, NC 28223-0001, USA}

\address[Alejandro Ram\'\i rez]{Facultad de Matem\'aticas\\
Pontificia Universidad Cat\'olica de Chile\\ Vicu\~na Mackenna 4860, Macul\\
Santiago 6904411, Chile\newline
Url: {\rm http://www.mat.puc.cl/\~\ \!\!\!\! aramirez}}

\bigskip
\email{benarous@cims.nyu.edu, smolchan@math.uncc.edu,aramirez@mat.puc.cl}

\begin{document}

\begin{abstract}

We describe a universal transition mechanism between annealed and
quenched regimes in the context of reaction-diffusion in random
media. We study the total population size for random walks which
branch and annihilate on  ${\mathbb Z}^d$, with time-independent
random rates. The random walks are independent,  continuous time,
rate $2d\kappa$, simple, symmetric, with $\kappa \ge 0$. A random
walk  at  $x\in{\mathbb Z}^d$,  binary branches at rate $v_+(x)$,
 and  annihilates at rate $v_-(x)$. The random
environment $w$ has coordinates
$w(x)=(v_-(x),v_+(x))$ which are i.i.d.
We identify a natural way to describe the annealed-Gaussian transition
mechanism under mild conditions on the rates.
Indeed, we introduce the exponents
$F_\theta(t):=\frac{H_1((1+\theta)t)-(1+\theta)H_1(t)}{\theta}$,
 and  assume that
$\frac{ F_{2\theta}(t)-F_\theta(t)}{\theta\log(\kappa t+e)}\to\infty$ for
$|\theta|>0$ small enough,
where
 $H_1(t):=\log \langle m(0,t)\rangle$ and $\langle m(0,t)\rangle$
denotes the average  of the expected value of the
number of particles $m(0,t,w)$ at time $t$ and an
environment of rates $w$, given that initially
there was only one particle at $0$.
 Then
  the empirical average of $m(x,t,w)$ over a box of
side $L(t)$ has different  behaviors: if  $ L(t)\ge e^{\frac{1}{d} F_\epsilon(t)}$ for some
$\epsilon >0$ and large enough $t$, a law of large numbers is satisfied;
if $ L(t)\ge e^{\frac{1}{d} F_\epsilon (2t)}$ for some
$\epsilon>0$ and large enough $t$,
a CLT is satisfied.
These statements are violated if the reversed inequalities
are satisfied for some negative $\epsilon$.
As corollaries, we obtain  more explicit statements
 under regularity conditions
on the tails of the random rates, including
examples  in the four universality classes
defined  in \cite{hkm}:  potentials which are
unbounded  of
Weibull type,    of double
exponential type,
 almost bounded,
 and
bounded  of Fr\'echet type. For them we also derive sharper results
in the non-annealed regime.

\end{abstract}

\renewcommand{\theenumi}{\alph{enumi}}
\renewcommand{\labelitemi}{}

    \makeatletter
    \def\thebibliography#1{\section*{References\@mkboth
      {REFERENCES}{REFERENCES}}\list
      {[\arabic{enumi}]}{\settowidth\labelwidth{[#1]}\leftmargin\labelwidth
        \advance\leftmargin\labelsep
        \usecounter{enumi}}
        \def\newblock{\hskip .11em plus .33em minus .07em}
        \sloppy\clubpenalty4000\widowpenalty4000
        \sfcode`\.=1000\relax}
    \makeatother

\newtheorem{proposition}{Proposition}
\newtheorem{theorem}{Theorem}
\newtheorem{remark}{Remark}
\newtheorem{model}{Model}
\newtheorem{lemma}{Lemma}
\newtheorem{corollary}{Corollary}

\def\sn{|\!|}
\maketitle

\vskip.4cm

\section{Introduction}

When studying the long time behavior of markovian dynamics in random
media, one is faced with an important distinction: quenched vs
averaged estimates. Should one work in the so-called quenched
regime, where the randomness of the medium is frozen, i.e. where the
dynamics are studied in one fixed random realization of the
medium/environment? Or should one work in the averaged regime where
both the randomness of the dynamics and of the medium are
considered, i.e. when one studies the dynamics in a given
realization but then also averages over the randomness of the
medium? \footnote{This regime is often called the annealed regime in
the mathematical literature. This is an entrenched but misleading
vocabulary since it is not the usual convention in physics}.

The true significance of this distinction "quenched vs averaged" is
important when these two regimes give different answers, which is
the case in many situations where the extreme values of the random
environment might play an important role. A good generic class of
examples where this distinction is significant is given by models of
reaction-diffusion in random media (\cite{gm}, \cite{gm2},
\cite{sznitman}).

There are two opposing views about this distinction. The  first
approach is to think that the relevant and important asymptotic
long-time estimates are the quenched ones. But, recognizing the
obvious fact that the quenched estimates  are the most intractable
ones mathematically, the averaged/annealed results are seen as a
welcome first approximation. The second and opposing view is that,
in many applications, the averaged asymptotic estimates are the only
relevant ones. The quenched results, though mathematically more
challenging, are not seen as useful or relevant. This second view is
naturally based on the idea that some mechanism must be at work,
which allows for averaging in the medium randomness.

Some years ago, the authors of this paper held the two opposing
opinions expressed above, based on their former collaborations with
different domains of applications (specifically physics of pollution
by underground waste storage  on the one hand and chemical kinetics
on the other). In the recent years we have built a common answer,
which we believe provides not only a natural resolution of the
scientific debate, but also introduces the idea that there exists in
fact a very rich universal transition between the averaged and
quenched results, which goes far beyond the reaction-diffusion
context. This transition also explains, in our view,the true
relevance of these regimes in various applications.

The key idea is the following: one should work with a fixed
realization of the medium but introduce a new parameter, say $L$,
which will be the scale of the spatial extent of the initial
distribution of the dynamics. Then, depending on the respective
sizes of the time scale $t$ and the space scale $L$ (when both t and
L diverge), one should see the transition we mentioned between the
quenched and averaged asymptotic results. More precisely one should
expect the following transition: for time scales $t(L)$ short enough
the averaged asymptotic results should be valid, whereas for very
long time scales $t(L)$ the quenched results should hold, and our
new intermediate asymptotics should emerge in between. The mechanism
of this transition is the following: for short time scales $t(L)$,
or equivalently for large space scales $L(t)$ for the initial
distribution,  the spatial ergodic theorem should ensure the needed
averaging mechanism in order to enforce the validity of
averaged/annealed results. For very large time scales no averaging
is possible and the extreme values of the environment play the
prominent role. The transition regimes are then easily understood,
they consist in regions of parameters $t$ and $L$ where one sees a
gradual emergence of the extreme values versus the average .

This scheme has been first established in the simplest possible
context, i.e i.i.d samples \cite{bbm}. More precisely sums of
exponential of i.i.d random variables are shown to exhibit this
transition.  In the context of reaction diffusion this could simply
be seen as reaction with no diffusion! In this simple context a full
transition is given from the Gaussian asymptotics to extreme value
theory: one sees in \cite{bbm} the gradual emergence of the
importance of extreme values of the i.i.d sample which gradually
destroys the validity of the Central Limit Theorem and then of the
Law of large Numbers by enforcing $\alpha$-stable fluctuations where
the exponent $\alpha$ decreases through the whole possible range,
i.e from $2$ to $0$. This mechanism is analogous to the phase
transition description of mean-field spin-glass equilibrium models
such as the Random Energy Model (\cite{d}, \cite{bkl}).

We then proceeded (\cite{bmr}) to one important case of
reaction-diffusion, i.e annihilation (or absorption) for random
walks in a random environment, more precisely random walks killed on
random obstacles, building on the work of \cite{sznitman}. There, we
studied the same transition mechanism for the natural quantity,
which is the probability of survival. Our picture was less precise
than in the i.id context, in the sense that even though we get the
proper scales for the intermediate regimes we cannot establish the
stable nature of the fluctuations, due to a lack of precise enough
understanding of the edge of the spectrum (for the generator of the
dynamics which is the discrete Dirichlet Laplacian on a random
domain of ${\mathbb Z}^d$).

In this paper we address a rather general case of reaction-diffusion
in random environments, i.e of a system of non-interacting
continuous-time Random Walks on the lattice ${\mathbb Z}^d$
branching and annihilating with stationary random rates (\cite{gm},
\cite{gm2}).

Let us first describe the random environment $\{w(x):x\in{\mathbb
Z}^d\}$, with $w(x):=(v_-(x),v_+(x))$ where
$v_+:=\{v_+(x):x\in{\mathbb Z}^d\}$ represents the branching rates
with $v_+(x)\in [0,\infty)$, while $v_-:=\{v_-(x):x\in{\mathbb
Z}^d\}$ represents the annihilation rates, with $v_-(x)\in
[0,\infty]$ so that we admit the possibility of  hard core
obstacles. We assume that the random variables $\{w(x):x\in{\mathbb
Z}^d\}$ are i.i.d and call their distribution $\mu$.

We consider the following dynamics in a fixed random realization of
this environment.
 We start with one particle at each site of the box $\Lambda_L$ of side $L$ in ${\mathbb Z}^d$.
 Each random walk moves
independently of the others according to a
 continuous time simple symmetric rate $2d\kappa$ dynamics
for some $\kappa\ge 0$ (we admit the possibility that $\kappa=0$,
i.e. no diffusion at all, as in \cite{bbm}). A random walk at a site
$x\in{\mathbb Z}^d$, branches at rate $v_+(x)$, disappearing and
producing two new independent offsprings,
 and annihilates at  rate $v_-(x)$ (note that
$v_-(x)=\infty$ means  the annihilation is instantaneous and certain
as in \cite{bmr}).

We study the asymptotic behavior, when both $t$ and $L$ go to
infinity, of the following observable of our system of random
walkers in random environment:
$$m_L(0,t,w):=\frac{1}{|\Lambda_L|}\sum_{x\in\Lambda_L}m(x,t,w)$$
Note that $m_L(0,t,w)$ is simply the (normalized) size of the total
population.

Our main result, theorem \ref{te1}, describes sharp conditions on
the scales of $t$ and $L$ for the validity of annealed asymptotics
of $m_L(0,t,w)$  and for these annealed asymptotics to cease to be
true. It also describes the scales for $t$ and $L$ where the
fluctuations about these annealed asymptotics are Gaussian and when
they cease to be so. We intend to give a full and complete picture
of the expected stable fluctuations, all the way to quenched
asymptotics, for a wide class of branching rates in a forthcoming
work.

The results we obtain here are general, in the sense that they are
valid for a large class of product distributions for the random
environment, under  mild conditions on the branching and
annihilation rates. They include  examples in the four universality
classes recently introduced by van den Hofstad, K\"onig and
M\"orters \cite{hkm} describing all the cases of the random
environment. For example,  potentials which are unbounded of Weibull
type, of double exponential type, almost bounded
 and bounded of Fr\'echet type.
Furthermore, our results include and generalize both theorems 2.1
and 2.2 of \cite{bbm},  parts $(ii)$ and $(iii)$ of theorem 1
 and part $(i)$ of theorem 2 of \cite{bmr}.

Let us now be more precise about our assumptions about the
distribution of the environment. In the recent paper \cite{hkm}, it
is shown that under regularity assumptions on the tail of the law of
the {\it effective potential} $v(0):=v_+(0)-v_-(0)$, exactly four
universality classes of environments can occur. Their assumptions
are formulated in terms of the {\it cumulant generating function},

\begin{equation}
\label{gen-cum}
H(t)=\log\left\langle e^{v(0)t}\right\rangle,\qquad\qquad t\ge 0,
\end{equation}
of the law of the effective potential $v(0)$, where for any function
of the environment $f$, we define $\langle f\rangle:=\int fd\mu$. Their basic
 assumption is that this function
is defined and finite for every $t\ge 0$. Then, under two regularity
assumptions on $H$ they show that four universality classes
can occur: $(1)$ a first class where $v$ is unbounded and has
``heavy tails'' at infinity, and which includes Weibull-type
tails; $(2)$ a second class of unbounded potentials with
``lighter'' tails which includes the double exponential law;
$(3)$ a class containing bounded and unbounded potentials;
$(4)$ a class of bounded potentials including those which
have Fr\'echet-type tails near their essential supremum,
and the degenerate case of random walks on hard core random obstacles
with $\mu[v(0)=-\infty]=p<1$.
In this paper we generalize parts $(ii)$ and $(iii)$ of theorem 1
 of \cite{bmr} describing the passage to an annealed and Gaussian
regime, to the previously described system of random walks
on the random environment $w$, under
mild conditions which include cases in these four universality classes.

Throughout this article the following will be assumed.

\smallskip
\noindent{\bf Assumption (E)}. {\it The law of the effective
potential is such that $\mu[v(0)=-\infty]<1$ and

$$
\left\langle e^{v(0)t}\right\rangle<\infty,\qquad t\ge 0.
$$}
\medskip

\noindent  Assumption {\bf (E)} ensures
 that $\mu$-a.s. the stochastic
process  of random walks on the
random  environment $w$
can be constructed on infinite volume (\cite{gm}), as a limit
of processes defined on finite boxes  corresponding to continuous
time Markov branching processes, as defined in Athreya and Ney
in \cite{an}. Furthermore, if $\zeta(t)$ denotes
the total number of random walks at time $t$
on a random environment $w$, given that
initially there was only a single one at site $x$,
and $E_x^w$ the expectation defined by its law,
 condition {\bf (E)} ensures the existence of the first moment $m(x,t,w)
=E_x^w[\zeta (t)]$.
This is the content of proposition \ref{non-explosion} of this paper.
This first moment will be the central object of our study.
We will see that assumption {\bf (E)} ensures the finiteness
for $t\ge 0$ of the annealed first
moments $\langle m(0,t)\rangle=\int m(0,t,w)d\mu$, with which we
will state our main assumptions.
We will need to
 define the
{\it growth functions} $\{H_1(t):t\ge 0\}$ by

\begin{equation}
\label{ge}
H_1(t):=\log\langle m(0,t)\rangle.
\end{equation}
Our main assumption will be formulated with the help of
a family of functions
$\{F_\theta:\theta\in{\mathbb R}\}$ which we call
the {\it intermittency  exponents}, defined
for every
$\theta\ne 0$ and $t\ge 0$  as,

\begin{equation}
\label{ef}
F_\theta(t):=
    \frac{H_1((1+\theta)t)-(1+\theta)H_1(t)}{\theta}.
\end{equation}

\smallskip

\noindent {\bf Assumption (MI)}.
{\it For all $|\theta|>0$ small enough,

$$
\lim_{t\to\infty}\frac{F_{2\theta}(t)-F_\theta(t)}
{\theta \log t}
=\infty.
$$
}
\medskip

\noindent As it will be shown in corollary \ref{cor-interm},
this assumption implies the occurrence
of the so called {\it intermittent} behavior of
the random field $w$ \cite{gm}.
It encompasses examples falling in the
four universality classes of \cite{hkm}.

We will show that it is possible to
formulate an assumption directly in terms of the cumulant
generating function $H$, which
is sufficient for {\bf (MI)} to be satisfied, and
includes the first class of \cite{hkm}.
 For this we
need to define the {\it cumulant  exponents}, parametrized
by $\theta\ne 0$, for $t\ge 0$ as,

\begin{equation}
\label{cum}
G_\theta (t):=\frac{H((1+\theta)t)-(1+\theta)H(t)}{\theta}.
\end{equation}
By Jensen's inequality it can be seen that $G_\theta (t)\ge 0$ whenever $0<|\theta|\le 1$.
 We will see that condition {\bf (MI)} is satisfied whenever the following happens.

\medskip

\noindent {\bf Assumption (SI)}.
{\it For all $|\theta|>0$ small enough,

$$
\lim_{t\to\infty}\frac{G_{2\theta}(t)-G_\theta(t)}{\theta  t}=\infty.
$$}

\noindent Condition {\bf (SI)} includes the first universality
class of \cite{hkm} and can be viewed as a strong intermittency
requirement.
It implies
 $H(t)/t\gg 1$. Hence, using the
bounds $e^{H(t)-2d\kappa t}\le\langle m(0,t)\rangle\le e^{H(t)}$
(see for example   theorem 3.1 of G\"artner and Molchanov \cite{gm}),
we see that if
 {\bf (SI)} is satisfied, we have the asymptotics

\begin{equation}
\label{gm}
\log\langle m(0,t)\rangle\sim H(t),
\end{equation}
which is much faster than the
a.s. one (see \cite{gm2}).

As already mentioned the interest of the results of the present
paper are their generality. Namely, their are valid only under the
assumptions {\bf (E)} and
 {\bf (MI)}.
Part $(i)$ of theorem \ref{te1}, states that if for
 some $\epsilon>0$
we have $L(t)\ge e^{\frac{1}{d} F_\epsilon(t)}$ eventually in $t$,
the law
of large numbers $\frac{m^L}{\langle m\rangle}\sim 1$ is satisfied
in probability:
hence we have the annealed  behavior $\log m^L(0,t)\sim \log\left\langle
m(0,t)\right\rangle$.
On the other hand, if for some $\epsilon>0$
we have $L(t)\le e^{\frac{1}{d} F_{-\epsilon}(t)}$ eventually
in $t$,
  in probability
$\frac{m^L}{\langle m\rangle}\ll 1$.
Part $(ii)$  says that if for some $\epsilon>0$
we have $L(t)\ge e^{\frac{1}{d} F_\epsilon(2t)}$ eventually in $t$,
then $\frac{m^L-(2L+1)^d\langle m\rangle}{Var_\mu m^L}$
 converges in distribution
to a centered normal law of unit variance ${\mathcal N}(0,1)$,
where $Var_\mu$ denotes the variance.
 Also, if for some $\epsilon>0$
we have $L(t)\le e^{\frac{1}{d} F_{-\epsilon}(2t)}$ eventually
in $t$,  in
probability $\frac{m^L-(2L+1)^d\langle m\rangle}{Var_\mu m^L}\ll 1$.
This discussion is summarized in the table below.

\begin{center}
\bigskip
\bigskip

\hskip0.2cm\vbox{\offinterlineskip
\def\strut{\vrule height 10.5pt
  depth 5.5pt width 0pt}
\halign{
\strut\vrule\quad\hfil # \hfil \quad
  &\vrule\qquad\hfill        $#$ \hfill\quad & \vrule\quad\hfill$#$\hfill \vrule  \cr
\noalign{\hrule}
Annealed behavior      &d \log L(t)\ge F_\epsilon (t) & m^L/\langle m
\rangle\sim 1
 \cr
Non-Annealed behavior   & d \log L(t)\le
F_{-\epsilon}(t)  & m^L/\langle m
\rangle\ll 1
\cr
Gaussian behavior   &  d \log L(t)\ge F_{\epsilon}(2t)
& \frac{m^L-(2L+1)^d\langle m\rangle}{Var_\mu m^L}\to {\mathcal N}(0,1)
  \cr
Non-Gaussian behavior   &  d \log L(t)\le F_{\epsilon}(2t)
& \frac{m^L-(2L+1)^d\langle m\rangle}{Var_\mu m^L
}\ll 1
\cr
\noalign{\hrule}
}
}
\end{center}

\centerline{{\ninebf Table 1}: \ninerm Large time asymptotic behavior
of the averaged first moments}

\bigskip

\noindent Under an additional regularity assumption on the
intermittency exponents (assumption {\bf (RI)} of
subsection 2.4) it will be shown in corollary \ref{corollaryr}
that there exist two constants
 $\gamma_1$ and $\gamma_2$, called {\it transition exponents} and a function
$J(t):[0,\infty)\to [0,\infty)$ with $J(t)\gg 1$, called
the {\it growth exponent},
 describing more explicitly  the transition of
theorem \ref{te1}.
 Indeed, in this case, a law of large numbers is satisfied
 when $d\log L(t)\ge \gamma J(t)$ for some $\gamma>\gamma_1$;
the CLT when $d\log L(t)\ge \gamma J(t)$ for some $\gamma>\gamma_2$.
Furthermore, if $d\log L(t)\le \gamma J(t)$ for some $\gamma<\gamma_1$, the
law of large numbers is not satisfied, while if
 $d\log L(t)\le \gamma J(t)$ for some $\gamma<\gamma_2$, the CLT
is not satisfied.
Propositions \ref{weibull}, \ref{double}, \ref{third} and
\ref{frechet} give the explicit value
of $\gamma_1$, $\gamma_2$ and $J$ in the case of unbounded
potentials with Weibull type tails, unbounded potentials
with double exponential type tails,
potentials in the third universality class of \cite{hkm} and
bounded potentials with Fr\'echet type tails including the
pure hard core case. Table 2 below summarizes those results.

\bigskip

\begin{center}
\bigskip
\bigskip

\hskip0.2cm\vbox{\offinterlineskip
\def\strut{\vrule height 10.5pt
  depth 5.5pt width 0pt}
\halign{
\strut\vrule\quad\hfil # \hfil \quad
  &\vrule\qquad\hfill        $#$ \hfill\quad
  &\vrule\qquad\hfill        $#$ \hfill\quad
  &\vrule\qquad\hfill        $#$ \hfill\quad
 & \vrule\quad\hfill$#$\hfill \vrule  \cr
\noalign{\hrule}
POTENTIAL   &\!\!\!\!\!\!\!-\log\mu[v(0)>x]\!\!\!   & J(t) & \gamma_1 & \gamma_2
 \cr
\noalign{\hrule}
Weibull   &  x^\rho, \rho>1      & H(t) & \frac{1}{1-\rho} & 2^{1-\gamma_1}\gamma_1
 \cr
Double exponential & e^{x/\rho}   & t  & \rho & 2\rho
\cr
Third class example&  e^{x^2}  & \frac{t}{2\sqrt{\log t}}     & 1  &2
  \cr
Pure hard core& --   &  c_2 t^{\frac{d}{d+2}} & \frac{2}{d+2}  & 2^{1-\gamma_1}\gamma_1
\cr
\noalign{\hrule}
}
}
\end{center}

\bigskip

\centerline{{\ninebf Table 2}: \ninerm Transition and growth exponents
in the four universality classes. In the pure hard }
\centerline{\ninerm core case
$\mu[v(0)=-\infty]=p$ and $c_2$ is a constant depending on $p$ and
$d$}
\medskip

\noindent Also, in theorem \ref{te2}, we obtain sharper upper bounds
for the order of magnitude of the averaged first moments in the
non-annealed regime for the examples treated in propositions \ref{weibull},
\ref{double}, \ref{third} and \ref{frechet}. This theorem, generalizes
{\bf Case 3} of theorem 2 of \cite{bmr}.

The special case in which assumption {\bf (SI)}
is satisfied expressed as
corollary \ref{corollary2}. This includes the
case $\kappa=0$, corresponding to sums of i.i.d. random
exponentials where condition {\bf (MI)}
reduces to,

\begin{equation}
\label{condition4}
\lim_{t\to\infty}\frac{G_{2\theta}(t)-G_\theta(t)}
{\theta}
=\infty,
\end{equation}
for $|\theta|>0$ large enough,
and $H_1(t)=H(t)$.
Corollary \ref{corollary2} is a result complementary to
theorem of \cite{bbm}, generalizing
theorems 2.1 and 2.2 of \cite{bbm} where Weibull and Fr\'echet type tails
are assumed on $v(0)$.

One of the main ingredients of the proof of theorem \ref{te1} is a
coarse graining technique, necessary to reduce the asymptotics of
the averaged first moments, to a sum of independent random
exponentials. This technique,
 was  introduced in \cite{bmr}, but here it faces the extra difficulty that the
terms of the sum defining the averaged probabilities are not
uniformly bounded with respect to the time variable $t$ or
the scale $L$ (whereas in \cite{bmr}, such a bound existed
having the value $1$). This requires more careful estimates
on these quantities, which are performed, via the
Feynman-Kac formula and spectral estimates. Once the reduction to a sum of
i. i. d. exponentials is achieved, an analysis based
on von Bahr-Esseen inequality finishes the proof (see
also \cite{bmr} and \cite{bbm}).

Besides this introduction, this paper has four other sections.
The main results are stated in section 2. We first introduce
in subsection 2.1 the main definitions. In subsection 2.2 we formulate
proposition \ref{non-explosion},
stating that a growth of the form $\limsup_{|x|\to\infty}
\frac{v_+(x)}{|x|\log |x|}=0$ is enough to ensure well defined
first moments for the total number of particles.
 When this proposition is combined
with proposition \ref{prop-gm} of
\cite{gm}, one concludes that  under the condition {\bf (E)},
the reaction-diffusion process on the lattice ${\mathbb Z}^d$
is such that the
total number of particles at any given time has a finite
first moment, for initial conditions with a finite total
number of particles. In particular, there is no explosion, and
the process is well defined.  We then state theorem \ref{te1}
in section 2.3. Corollary \ref{corollary2}, under
the assumption
{\bf (SI)} is stated and proved next. The applications
of theorem \ref{te1} are given in subsection 2.4.
First,  the regularity condition
 {\bf (RI)},
 is introduced.
This is applied
to  the case of unbounded effective
potentials with Weibull type tails, through
proposition  \ref{weibull}, using  the
Kasahara exponential Tauberian theorem \cite{bgt}.
Next, corollary \ref{corollaryr} is applied to
the case of unbounded potentials with double exponentially
decaying type tails, through proposition \ref{double}.
Then,  we treat the case of potentials
falling in the third universality class (almost bounded) of \cite{hkm}
through proposition \ref{third}. We end subsection 2.4
 considering the case of bounded potentials with
 tails near their upper-bound
which are of the Fr\'echet type (proposition \ref{frechet}).
In subsection 2.5, we state theorem \ref{te2}, which improves the
upper bounds describing the order of magnitude of the empirical
average for the examples considered in subsection 2.4.
The proof of proposition 1, is the content of the third section.
In section 4, the truncated first moments are introduced.
These
are the first moments of a reaction-diffusion process
defined on a finite box, with Dirichlet boundary conditions.
They are then used to approximate
some important quantities related to the averaged first
moments.
Then,  several important estimates for the moments and correlations
of the first moments are derived. The proof of theorems \ref{te1} and
\ref{te2} are given in section 5.
In subsection 5.1, the partition analysis method of \cite{bmr}
is recalled. This and together with the estimates
of section 4, and the von Bahr-Esseen inequality, is subsequently
applied to prove
theorem \ref{te1}. The paper finishes with subsection 5.7, where
theorem \ref{te2} is proved.

\vskip.7cm

\section{Notation and results}

Here we will state the results of this paper. After introducing
most of the notation and giving
the main definitions in the first subsection, in subsection 2.2
we state proposition  \ref{assumet}, which ensures that a.s. there is no
explosion for the reaction-diffusion process under assumption
{\bf (E)}. Then, the principal result of this paper, theorem
\ref{te1}, is stated in subsection 2.3, together with corollary
\ref{corollary2}.
In subsection 2.4, we state corollary \ref{corollaryr},  giving the form of
theorem \ref{te1}, under certain regularity assumptions. Here we will
consider applications  of this results to several specific
examples of distributions of the effective potential. We end the
presentation of our results with subsection 2.5, where we state
theorem \ref{te2}.

\vskip.4cm

\subsection{Definition of the reaction-diffusion process}
We begin defining a reaction-diffusion model corresponding to a set of
random walks on the lattice ${\mathbb Z}^d$ branching
 and annihilating at rates depending on the their position. Consider the set of natural numbers ${\mathbb N}$ endowed with the
discrete topology.
Define the set $\Omega:={\mathbb N}^{{\mathbb Z}^d}$  representing the
possible configuration of particles on the lattice. In this paper we will
be interested only on the subset of configurations $\Omega_0\subset\Omega$
characterized by the property that $\{x:\eta(x)> 0\}$ has
finite cardinality whenever $\eta\in\Omega_0$.
 Let  $v_+:=\{v_+(x):x\in{\mathbb Z}^d\}$
and $v_-:=\{v_-(x):x\in{\mathbb Z}^d\}$,
where $v_+(x)\in [0,\infty)$ and $v_-(x)\in [0,\infty]$.
 Here $v_+(x)$ and $v_-(x)$ represent the rate at which
particles branch and annihilate at site $x$, respectively.
Note that we admit the value $\infty$ for the annihilation rate:
this represents a hard core
obstacle, where particles are instantly annihilated.
 Call
an ordered pair $w:=(v_-,v_+)\in W$, with
coordinates $w(x)=(v_-(x),v_+(x))$, a
{\it field configuration}, where $W:=([0,\infty)\times [0,\infty])^{{\mathbb Z}^d}$.
We will denote the set of {\it hard core obstacle sites}
by ${\mathcal G}(w):=\{x\in{\mathbb Z}^d:v_-(x)=\infty\}$.
Given $r\in [0,\infty)$ and $x\in{\mathbb Z}^d$  we will call
$\Lambda (x,r):=\{y\in{\mathbb Z}:\sn x-y\sn\le r\}$ the ball of radius
$r$ centered at $x$ under the norm $\sn x\sn :=\sup_{1\le i\le d}|x_i|$, where
$x_i$ are the coordinates of $x$. We will furthermore use the
notation $\Lambda_r$ in place of $\Lambda(0,r)$.
In this subsection we will construct a stochastic process as the limit
as $L\to\infty$ of processes defined on  the boxes
 $\Lambda_L$. Throughout the sequel,
given a subset $U\subset{\mathbb Z}^d$, we will denote by
$U^c$ the complement of $U$ and $\delta U:=\{x\in U^c:\inf_{y\in U}|x-y|=1\}$,
where $|\cdot |$ denotes the Euclidean distance.
So, for each finite subset $U\subset{\mathbb Z}$, we want to
 consider a process
 with state space
$\Omega^w_0:=\{\eta\in\Omega_0:\eta(x)=0\ ,x\in{\mathcal G}(w)\}$ defined formally by the infinitesimal generator,

\begin{eqnarray}
\nonumber
L_U f(\eta)&:=&\sum_{x\in U\cap{\mathcal G}(w)^c}
\sum_{y\in{\mathcal G}(w)^c:||x-y||=1}
  \kappa \eta(x) (f(\eta^{x,y})-f(\eta))\\
\nonumber
&+&\sum_{x\in U\cap{\mathcal G}(w)^c}
\sum_{y\in{\mathcal G}(w):||x-y||=1}
  \kappa \eta(x) (f(\eta^{x,-})-f(\eta))\\
\nonumber
&+&\sum_{x\in U\cap{\mathcal G}(w)^c} v_+(x)\eta(x)(f(\eta^{x,+})-f(\eta))\\
\label{rd1}
&+&\sum_{x\in {\mathcal G}(w)^c} v_-(x)\eta(x)(f(\eta^{x,-})-f(\eta)),
\end{eqnarray}
acting on an appropriate dense subset ${\mathcal D}_U$ of the
space of real bounded functions $B(\Omega^w_0)$,
defined on $\Omega^w_0$, endowed with the uniform norm.
In the above expression,
$\kappa\ge 0$, $\eta^{x,y}$ is the configuration where a particle
from site $x$ has jumped to site $y$ so that
$\eta^{x,y}(z)=\eta(z)$ if $z\ne x,y$, $\eta^{x,y}(x)=\eta(x)-1$, and
$\eta^{x,y}(y)=\eta(y)+1$;
$\eta^{x,+}$ is the configuration where there is an extra particle
at site $x$ and $\eta^{x,-}$ the configuration where one particle at site $x$
has disappeared. It is a well known fact that it is possible
to construct a strong Markov process, denoted by $\eta^U:=\{\eta^U (t):t\ge 0\}$, corresponding to an infinitesimal generator of the form
(\ref{rd1}), and taking values on the
 Skorokhod space ${\mathcal S}:=D([0,\infty);\Omega^w_0)$. In fact, such
a process falls in the category called {\it $|U|$-dimensional
continuous time Markov  branching process}
by Athreya and Ney (see chapter V, sections 7.1-7.2 of Athreya-Ney \cite{an}).
Furthermore, it can be shown that a.s. the expected value of
each coordinate of such a process $\{\eta^U(t):t\ge 0\}$, is finite,
ensuring that there cannot be infinitely many particles produced in a finite
 time (see section 7.1 of \cite{an}).
Let us now call
${\mathcal P}(\Omega^w_0)$ the set of probability measures defined on
$\Omega^w_0$ endowed with the Borel $\sigma$-algebra associated to
the subspace topology of $\Omega^w_0$ as a subset of
 $\Omega$ with   the product topology.
Then, for
each field configuration $w\in W$
and  probability measure $\nu\in{\mathcal P}(\Omega^w_0)$
denote by $P_{\nu}^{U,w}$ the law of the process $\{\eta^U(t):t\ge 0\}$
defined on ${\mathcal S}$ endowed with its Borel $\sigma$-field
${\mathcal B}({\mathcal S})$.
We will call this process the {\it reaction-diffusion process
on $U$} with field $w$ and initial condition $\nu$.
In the particular case in which $U=\Lambda_n$ we will use
the obvious notations $L_n$ and $P_{\nu}^{n,w}$.
Furthermore, we will call the process on $\Lambda_n$, the
{\it reaction-diffusion process
at scale $n$}.
Now, note that using the natural
coupling \cite{liggett} and Kolmogorov's extension theorem,
 it is possible to define for each field configuration
$w\in W$ and initial condition $\nu\in{\mathcal P}(\Omega^w_0)$
a probability measure $Q_\nu^w$ on the product space
${\mathcal S}^{\mathbb N}$, endowed with its Borel $\sigma$-algebra
 induced by the product topology, in such a way that
if $\eta^n\in{\mathcal S}$ denotes the $n$-th coordinate
of an element $\eta\in{\mathcal S}^{\mathbb N}$, $\eta^n(t)\in\Omega^w_0$
its value at time $t\ge 0$ and $\eta^n(t,x)\in{\mathbb N}$ the value at
time $t$ of the $x$-coordinate of $\eta^n(t)$, then,

\begin{itemize}
\item[(i)] for every $A\in{\mathcal B}({\mathcal S})$ and $n\ge 1$,

$$Q_\nu^w[\eta^n\in A]=P_{\nu}^{n,w} [A].$$
In particular, for every $B\in{\mathcal B}(\Omega_0)$ we have that,
$Q_\nu^w[\eta^n(0)\in B]=\nu [B].$

\item[(ii)] for every $n\ge 1$,

$$Q_\nu^w[\eta^{n+1}(t)\ge\eta^n(t)]=1.$$
\item[(iii)] for each $n\ge 1$ define the first exit time from the
box $\Lambda_n$ as

$$T_n:=\inf\{t\ge 0:\sup_{x\in\Lambda_n^c}\eta^n (x,t)>0\}.$$
Then, for every $n\ge 1$,

$$Q_\nu^w[\eta^n(t)=\eta^{n+1}(t), T_n>t]=Q_\nu^w[T_n>t].$$

\end{itemize}
Let us now remark that due to property $(ii)$, for every $t\ge 0$
and $x\in{\mathbb Z}^d$ the limit,

$$\eta(t,x):=\lim_{n\to\infty}\eta^n(t,x),$$
exists, possibly taking the value $\infty$. Define $\eta(t):=\{\eta(t,x):
x\in{\mathbb Z}^d\}$. We denote the stochastic process
$\{\eta(t):t\ge 0\}$, taking values on the space
$\bar{\mathbb N}^{{\mathbb Z}^d}$, where $\bar{\mathbb N}$ is the
 Alexandrov compactification
of the natural numbers, and distributed according to the
measure $Q_\nu^w$, the  {\it reaction-diffusion process with field $w$ and
initial law $\nu$}. We will denote by $E_\nu^w$ the
corresponding expectation.
In addition, for each $t\ge 0$, we define the {\it total number of particles
at time $t$} by

$$\zeta(t):=\sum_{x\in{\mathbb Z}^d}\eta(t,x).$$
Also, whenever it is true that,

$$Q^w_\nu\left[ \eta(x,t)<\infty,\ \forall\ x\in{\mathbb Z}^d,\ t\ge 0\right]
=1,$$
we will say that with probability one there is {\it no explosion}.
In the sequel we define for each $x\in{\mathbb Z}^d$ the probability
measure $\delta_x$ on $(\Omega^w_0,{\mathcal B})$ which assigns probability $1$
to configurations with one particle at site $x$ and none elsewhere.
We will be interested in initial configurations where $\nu=\delta_x$.
In such a case we will use the notation $P^w_x$ instead of $P^w_{\delta _x}$
and $E^w_x$ for the corresponding expectation. In the case
where $x\in{\mathcal G}(w)$, we adopt the convention that
$P^w_x$ is the probability measure which has a unique
atom at the configuration $\eta\equiv 0$ ($\eta(x)=0$ for every
$x\in{\mathbb Z}^d$).

Let us denote by ${\mathcal P}(W)$ the set of probability measures
defined on the space
$W=([0,\infty)\times [0,\infty])^{{\mathbb Z}^d}$ endowed with
its natural $\sigma$-algebra.
In the sequel we will take fields $v_+$ and $v_-$ which are random,
assigning  a probability measure $\mu\in{\mathcal P}(W)$
 in such a way that the field configuration
 $\{w(x):x\in{\mathbb Z}^d\}$
has independent coordinates with respect to $\mu$.
 Furthermore, we will use the notation $\langle\cdot\rangle$ to
denote expectation with respect to this law and $Var_\mu(\cdot )$
variance.
Now, let us define
 the {\it quenched first moment} on
${\mathbb Z}^d$ of the total number of particles at time
$t$ starting from site $x$ as,
$m(x,t,w):=E^w_x\left[\zeta (t)\right]$,
and the {\it annealed first moment} on ${\mathbb Z}^d$
of the total number of particles at time $t$ starting from site $x$ as,
$\langle m(x,t)\rangle:=\int m(x,t,w)d\mu$.  Furthermore, we call the
sets $\{m(x,t,w):x\in{\mathbb Z}^d\}$ and
$\{\langle m(x,t)\rangle:x\in{\mathbb Z}^d\}$,  the
{\it fields of quenched first moments} and the {\it field of annealed
first moments}, respectively. Depending on the context, we
might write $m$ or $m(x,t)$ in place of $m(x,t,w)$, dropping
the dependence on the field configuration $w$, and $\langle m\rangle$
instead of $\langle m(x,t)\rangle$.

The quantity which will give us a transition mechanism between
of the  quenched first moments is the
{\it averaged first moment} at scale $L$ and time $t$ defined
for a reaction-diffusion process starting from site $x$ as,

$$ m^L(x,t,w):=\frac{1}{|\Lambda(x,L)|}\sum_{y\in\Lambda(x,L)}m(y,t,w).$$

\smallskip

\subsection{Results giving conditions for no explosion}
Here we will give a criteria on the field configuration $w$,
stated as proposition \ref{non-explosion}, which
ensures that there is no explosion in the reaction-diffusion process
with field $w$.

\medskip
\begin{proposition}
\label{non-explosion}
 Consider the reaction diffusion process
with field $w$ and initial law $\nu$. Assume that,

\begin{equation}
\label{assumet}
\limsup_{|x|\to\infty}\frac{v_+(x)}{|x|\log |x|}=0.
\end{equation}
Then for every $t\ge 0$ we have that

$$E^w_\nu[\zeta(t)]<\infty.$$
Hence, there is no explosion.
\end{proposition}

\medskip

\noindent
We state below with the name of proposition \ref{prop-gm}, a
 result of G\"artner-Molchanov \cite{gm} (lemma 2.5) giving a sufficient
condition on the law $\mu$ in order that the first moment of the
total number of particles $\zeta(t)$ at time $t$, exists $\mu$-a.s.,
and hence that there is  no explosion.
 Given a field configuration $w=(v_+,v_-)$, we
 now need to introduce the
{\it effective potential} $\{v(x):x\in {\mathbb Z}^d\}$, defined by
$v(x):=v_+(x)-v_-(x)$. Furthermore, set $\log_+ x=\log x$ if $x>e$ and
$\log_+ x=1$ otherwise, while define the positive part $x^+:=\max (0,x)$.

\medskip
\begin{proposition}
\label{prop-gm} Consider the reaction diffusion process
with field $w$ and initial law $\nu$. Assume that the
field configuration $w$ is distributed according to a product
probability measure $\mu\in {\mathcal P}(W)$. Suppose
that

\begin{equation}
\label{eee}
\left\langle\left(\frac{v^+(0)}{\log_+ v(0)}\right)^d\right\rangle<\infty
\end{equation}
Then $\mu$-a.s. it is true that,

$$\limsup_{|x|\to\infty}\frac{v_+(x)}{|x|\log |x|}= 0,$$
and therefore, $\mu$-a.s. there is no explosion for the reaction-diffusion process with field $w$ and arbitrary initial law
in ${\mathcal P}(\Omega_0)$.
\end{proposition}

\noindent Note that the last statement of proposition \ref{prop-gm}, follows
from proposition \ref{non-explosion}. Furthermore, assumption {\bf (E)} is enough
for (\ref{eee}) to be satisfied, and hence to ensure no explosion.
\smallskip

\subsection{The Gaussian-Annealed transition results} Here we will
state the main result of this paper, which shows how under
different growth of scales, the averaged first moment has an
asymptotic behavior where a law of large numbers is satisfied,
 and a central
limit theorem can describe the fluctuations around this law
of large numbers. We will assume condition {\bf (E)}, ensuring
the existence of the annealed first moments.
We will also need to
consider the {\it growth functions} $\{H_1(t)=\log\langle m(0,t)\rangle:t\ge 0\}$
already defined in display (\ref{ge}) of the introduction
and the   {\it intermittency  exponents} $\{F_\theta:\theta\in{\mathbb R}\}$, defined in display (\ref{ef}).
Let us now state the main result of this paper.

\medskip

\begin{theorem}
\label{te1} Consider a reaction-diffusion process
with initial law $\delta_0$ and
field $w=(v_+,v_-)$ distributed according to a product measure
$\mu\in{\mathcal P}(W)$.
 Consider the intermittency  exponents $\{F_\theta:\theta\in{\mathbb R}\}$
defined in display (\ref{ef}) and the growth functions
 $\{H_1(t):t\ge 0\}$ defined in display (\ref{ge}). Assume that conditions
 {\bf (E)} and {\bf (MI)} are
 satisfied.
Then the following statements are true,

\begin{itemize}
\item[(i)] {\bf Law of large numbers}. Assume that there is an $\epsilon>0$
such that  eventually in $t$,
$L(t)\ge \exp\left\{\frac{1}{d} F_{\epsilon}(t)\right\}$.
Then in $\mu$-probability we have

\begin{equation}
\label{int1}
\frac{ m^L(0,t,w)}{\langle m(0,t)\rangle}\sim 1,
\end{equation}
as $t\to\infty$.
Furthermore, assume that there is an $\epsilon>0$ such that eventually
in  $t$,
$L(t)\le \exp\left\{\frac{1}{d} F_{-\epsilon} (t)\right\}$.
Then in $\mu$-probability we have

\begin{equation}
\label{int2}
\frac{ m^L(0,t,w)}{\langle m(0,t)\rangle}\ll 1,
\end{equation}
as $t\to\infty$.

\item[(ii)] {\bf Central limit theorem}. Assume that there is
an $\epsilon>0$ such that eventually in $t$,
$L(t)\ge \exp\left\{\frac{1}{d} F_{\epsilon} (2t)\right\}$. Then,

\begin{equation}
\label{int3}
\lim_{t\to\infty}\frac{ m^L(0,t,w)-\langle m(0,t)\rangle}{\sqrt{Var_\mu (m(0,t))}}={\mathcal N}(0,1),
\end{equation}
where ${\mathcal N}(0,1)$ is a centered normalized normal law and
the convergence is in the sense of distributions. Furthermore,
assume that there is
an $\epsilon>0$ such that eventually in $t$,
$L(t)\le \exp\left\{\frac{1}{d} F_{-\epsilon} (2t)\right\}$. Then
 in $\mu$-probability we have that,

\begin{equation}
\label{int4}
\frac{ m^L(0,t,w)-\langle m(0,t)\rangle}{\sqrt{Var_\mu(m(0,t))}}\ll 1.
\end{equation}
\end{itemize}
\end{theorem}
\smallskip

\noindent To provide a better insight on the
meaning of theorem \ref{te1},  we will see the form
that it takes under the stronger assumption {\bf (SI)},
which includes the first universality class of \cite{hkm}.
This will be formulated as a corollary,
which in the case $\kappa=0$ generalizes theorem 2.1 and 2.2
of \cite{bbm} to
include distributions $\mu$ of the field, which
not necessarily have regularly varying log-tails.

\begin{corollary}
\label{corollary2}
Consider a reaction-diffusion process
with initial law $\delta_0$ and
field $w=(v_+,v_-)$ distributed according to a product measure
$\mu\in{\mathcal P}(W)$.
 Consider the cumulant intermittency  exponents $\{G_\theta:\theta\in{\mathbb R}\}$
defined in display (\ref{cum}) and the cumulant generating
function
 $\{H(t):t\ge 0\}$ defined in display (\ref{gen-cum}). Assume that conditions
 {\bf (E)} and {\bf (SI)} are satisfied.
Then the following statements are true,

\begin{itemize}
\item[(i)] {\bf Law of large numbers}. Assume that there is an $\epsilon>0$
such that  eventually in $t$,
$L(t)\ge \exp\left\{\frac{1}{d} G_{\epsilon}(t)\right\}$.
Then in $\mu$-probability we have

\begin{equation}
\nonumber
\frac{ m^L(0,t,w)}{\langle m(0,t)\rangle}\sim 1,
\end{equation}
as $t\to\infty$. In particular, $\frac{ \log m^L(0,t,w)}{H(t)}\sim 1$.
Furthermore, assume that there is an $\epsilon>0$ such that eventually
in  $t$,
$L(t)\le \exp\left\{\frac{1}{d} G_{-\epsilon} (t)\right\}$.
Then in $\mu$-probability we have

\begin{equation}
\nonumber
\frac{ m^L(0,t,w)}{\langle m(0,t)\rangle}\ll 1,
\end{equation}
as $t\to\infty$.

\item[(ii)] {\bf Central limit theorem}. Assume that there is
an $\epsilon>0$ such that eventually in $t$,
$L(t)\ge \exp\left\{\frac{1}{d} G_{\epsilon} (2t)\right\}$. Then,

\begin{equation}
\nonumber
\lim_{t\to\infty}\frac{ m^L(0,t,w)-\langle m(0,t)\rangle}{\sqrt{Var_\mu (m(0,t))}}={\mathcal N}(0,1),
\end{equation}
where the convergence is in the sense of distributions. Furthermore,
assume that there is
an $\epsilon>0$ such that eventually in $t$,
$L(t)\le \exp\left\{\frac{1}{d} G_{-\epsilon} (2t)\right\}$. Then
 in $\mu$-probability we have that,

\begin{equation}
\nonumber
\frac{ m^L(0,t,w)-\langle m(0,t)\rangle}{\sqrt{Var_\mu(m(0,t))}}\ll 1.
\end{equation}
\end{itemize}
\end{corollary}

\noindent The proof of  corollary \ref{corollary2} in the case $\kappa =0$
follows from the observation that $G_\theta(t)=F_\theta(t)$ in this
case.
 The case $\kappa>0$ is a direct
consequence of the fact that $e^{H(t)-2d\kappa t}
\le \langle m(0,t)\rangle\le e^{H(t)}$, stated in
theorem 3.1 of \cite{gm}, and the observation that
$G_\theta (t)\ge 0$ for $\theta>0$ and $G_\theta (t)\le 0$
for $\theta<0$, which follows from Jensen's inequality.
Now, the following proposition, which will be proved in
section \ref{mce}, shows that the condition,

\begin{equation}
\label{condition3}
\lim_{t\to\infty}tH''(t)=\infty,
\end{equation}
is sufficient for assumption {\bf (SI)} to be true. This
condition implies a kind of domination of the branching
over the annihilation.
\smallskip
\begin{proposition}
\label{srel2}
 Consider the cumulant   exponents
$\{G_\theta:\theta\in{\mathbb R}\}$.
Assume that condition (\ref{condition3}) is satisfied.
Then, condition {\bf (SI)} is satisfied. Furthermore:
$(i)$ for every $\theta\ne 0$,
there is a $t_0\ge 0$, such that
 the function $G_\theta (t)$ is monotone in $t$, for $t\ge t_0$;
$(ii)$ there is a $t_1\ge 0$,
 such that the function $G_\theta (t)$ is monotone in $\theta$
for $t\ge t_1$.
\end{proposition}
\smallskip

\noindent Condition (\ref{condition3}) implies that
 the  branching dominates the annihilation, in the sense
that $H(t)/t\to\infty$ as $t\to\infty$, which implies
that the essential supremum of the
 random variable $v(0)$ is infinite
(see \cite{gm}).

\smallskip
\subsection{Regularity assumptions on the intermittency exponents}.
For the purpose of applications,
it will be important to identify cases where the assumptions in theorem
\ref{te1} can be formulated in a more explicit way.
As it will be shown, the following assumption on the intermittency exponents, turns out to fall
in one of these situations.

\medskip

\noindent {\bf Assumption (RI)}. {\it The intermittency exponents
$\{F_\theta(t):\theta\ne {\mathbb R}\}$ satisfy the mild
intermittency condition {\bf (MI)}. In addition, there
exist two
increasing functions $f_1, f_2:{\mathbb R}\backslash\{0\}\to {\mathbb R}$ and
a function $J(t):[0,\infty)\to [0,\infty)$, such that
for $\theta\ne 0$ small enough,

\noindent $(i)$ $
F_\theta(t)\sim f_1(\theta) J(t),
$ and $
F_\theta(2t)\sim f_2(\theta) J(t).
$

\noindent $(ii)$ There exists two constants $\gamma_1$ and $\gamma_2$,
such that $\lim_{\theta\to 0} f_1(\theta)=\gamma_1$ and
$\lim_{\theta\to 0} f_2(\theta)=\gamma_2$.

}
\medskip Throughout the sequel, the constants
$\gamma_1$ and $\gamma_2$ will be called {\it transition exponents}
and the function $J$ {\it growth exponent}. It will be shown that there exist several important
cases of random fields which fall in this category.
Furthermore, the following corollary of theorem \ref{te1} shows the convenience
of assumption {\bf (RI)}.

\smallskip

\begin{corollary}
\label{corollaryr}
Consider a reaction-diffusion process
with initial law $\delta_0$ and
field $w=(v_+,v_-)$ distributed according to a product measure
$\mu\in{\mathcal P}(W)$.
Suppose that assumptions {\bf (E)} and {\bf (RI)} are
 satisfied with transition exponents $\gamma_1$ and $\gamma_2$ and
growth function $J$.
 Then the following statements are true,

\begin{itemize}
\item[(i)] {\bf Law of large numbers}.
Assume that there is a
$\gamma>\gamma_1$
such that eventually in  $t$,
$\log L(t)\ge \frac{1}{d} \gamma J(t)$.
Then in $\mu$-probability we have

\begin{equation}
\nonumber
\frac{ m^L(0,t,w)}{\langle m(0,t)\rangle}\sim 1,
\end{equation}
as $t\to\infty$.
Furthermore, assume that there is a $0<\gamma<\gamma_1$
 such that eventually in  $t$,
$\log L(t)\le \frac{1}{d}\gamma J (t)$.
Then in $\mu$-probability we have

\begin{equation}
\nonumber
\frac{ m^L(0,t,w)}{\langle m(0,t)\rangle}\ll 1,
\end{equation}
as $t\to\infty$.

\item[(ii)] {\bf Central limit theorem}.
 Assume that there is
a $\gamma>\gamma_2$ such that eventually in $t$,
$\log L(t)\ge \frac{1}{d}\gamma J (t)$. Then,

\begin{equation}
\nonumber
\lim_{t\to\infty}\frac{ m^L(0,t,w)-\langle m(0,t)\rangle}{\sqrt{Var_\mu (m(0,t))}}={\mathcal N}(0,1),
\end{equation}
where
the convergence is in the sense of distributions. Furthermore,
assume that there is
a $0<\gamma<\gamma_2$ such that eventually in $t$,
$\log L(t)\le \frac{1}{d}\gamma J (t)$. Then
 in $\mu$-probability we have that,

\begin{equation}
\nonumber
\frac{ m^L(0,t,w)-\langle m(0,t)\rangle}{\sqrt{Var_\mu(m(0,t))}}\ll 1.
\end{equation}
\end{itemize}
\end{corollary}
\smallskip

\noindent In the next subsection we will apply corollary \ref{corollaryr}
to four situations each one falling in one of the universality
classes described be van den Hofstad, K\"onig and M\"orters
in \cite{hkm}. These classes encompass all possible situations
under three conditions. The first is condition {\bf (E)}, ensuring the existence
of the positive moments defining the cumulant generating functions
(\ref{gen-cum}).  The second and
third condition avoids different qualitative
behaviors of the potential at different scales. Le us formulate
next the second condition of \cite{hkm}.

\medskip

\noindent {\bf Assumption (H)}. {\it The function $H(t)/t$ is in the
de Haan class.}

\medskip

\noindent A function $f$ is said to be in the de Haan class if
for some regularly varying function $g:(0,\infty)\to {\mathbb R}$, we
have that $(f(\lambda t)-f(t))/g(t)$ converges to a limit
different from $0$ as $t\to\infty$, for
$\lambda>0$.
Let us recall that a function
$h$ is regularly varying at infinity with index $\rho$,
if for any $a>0$ we have $\lim_{x\to\infty} h(ax)/h(x)=a^\rho$.
This property will be stated as  $h\in R_\rho$.
 Whenever {\bf (H)} is satisfied, then
$H(t)$ is regularly varying with index $\gamma\ge 0$.
Furthermore, in proposition 1.1 of \cite{hkm},
it is proven that under assumption {\bf (H)} there exist
a function $\hat H:(0,\infty)\to {\mathbb R}$ and a continuous
function $k(t):(0,\infty)\to (0,\infty)$ such that,

\begin{equation}
\label{p11}
\lim_{t\to\infty}\frac{H(ty)-yH(t)}{k(t)}=\hat H(y)\ne 0,
\end{equation}
for $y\in (0,1)\cup (1,\infty)$. It is also shown that
$k(t)$ is regularly varying of index $\gamma$. We can now
recall the third assumption of \cite{hkm}.

\medskip

\noindent {\bf Assumption (K)}. {\it The limit $k^*=
\lim_{t\to\infty}\frac{k(t)}{t}$ exists in $[0,\infty]$.}

\medskip

\noindent Under assumptions {\bf (E)}, {\bf (H)} and
{\bf (K)}, the four universality classes defined in
\cite{hkm} are:

\begin{itemize}
\item[(1)] $\gamma>1$, or $\gamma=1$ and $k^*=\infty$.

\item[(2)] $\gamma=1$ and $k^*\in (0,\infty)$.

\item[(3)] $\gamma=1$ and $k^*=0$.

\item[(4)] $\gamma<1$.

\end{itemize}

\smallskip

\noindent In what follows we exhibit examples in each one
of these classes satisfying assumption {\bf (RI)} so that corollary \ref{corollaryr} can be
applied, and the transition exponents can be explicitly written.
In the sequel, following \cite{hkm}, we will call the third class
the class of {\it almost bounded potentials}.
\smallskip

\subsubsection{Unbounded potentials with Weibull type
 tails} Our first application of corollary \ref{corollary2}
will be to an example falling in the first universality class
of \cite{hkm}. We will assume that
the essential supremum of $v(0)$ is $\infty$, and the tails
at $\infty$ of $v(0)$ follow a Weibull-type law,
 $\mu[v(0)>x]
=\exp\{-h(x)\}$ for $x>0$, with $h\in R_{\rho}$ for some
$1<\rho<\infty$. In the terminology
of \cite{bbm} in the context of i.i.d. random exponentials
($\kappa=0$ in our situation), this is called {\it Case} {\bf B}.
The following proposition shows that assumption {\bf (RI)} is
satisfied in this situation, and hence corollary \ref{corollaryr},
which generalizes theorems 2.1 and 2.2 of \cite{bbm} in case
{\bf B} form $\kappa=0$ to $\kappa>0$.

\medskip

\begin{proposition}
\label{weibull}
 Suppose that the essential supremum
of the effective potential $v(0)$ is $\infty$ with
Weibull-type tails $\mu[v(0)>x]
=\exp\{-h(x)\}$ for $x>0$, and $h\in R_{\rho}$ for some
$1<\rho<\infty$. Then assumption {\bf (RI)} is satisfied with
transition exponents

$$\gamma_1=\frac{1}{\rho-1},\qquad\qquad
\gamma_2= 2^{\frac{\rho}{\rho-1}}\frac{1}{\rho-1},$$
and growth exponent,

$$J(t)=H(t).$$
\end{proposition}

\medskip

\noindent Let us now prove  proposition \ref{weibull}.
Note that under the conditions on the tail of the law
of $v(0)$, the cumulant generating function
$H(t)$ of $v(0)$ is well defined, smooth, non-decreasing and tends to
infinity as $t\to\infty$.  Furthermore, by the Kasahara exponential
Tauberian theorem (see Bingham, Goldie and Teugels
\cite{bgt}, Theorem 4.12.7), we
know that $H\in R_{\rho '}$, where the index $\rho'$ is defined
by the equation,

$$\frac{1}{\rho}+\frac{1}{\rho'}=1.$$
>From this observation it is easy to check that
assumption {\bf (SI)} of corollary \ref{corollary2} is
 satisfied, and that
 for
every $\epsilon\ne 0$, the cumulant growth exponents  $G_\epsilon (t)$
and $G_\epsilon (2t)$  satisfy

$$
\frac{G_\epsilon(t)}{H(t)}
\sim\frac{(1+\epsilon)^{\rho'}-(1+\epsilon)}\epsilon$$
and
$$
\frac{G_\epsilon(2t)}{H(t)}
\sim 2^{\rho'}\frac{(1+\epsilon)^{\rho'}-(1+\epsilon)}\epsilon.$$
Now note that $\frac{(1+\epsilon)^{\rho'}-(1+\epsilon)}\epsilon$ is
increasing in $\epsilon$ and converges to $\gamma_1=\frac{1}{\rho-1}$
as
$\epsilon\to 0$, and similarly
 $2^{\rho'}\frac{(1+\epsilon)^{\rho'}-(1+\epsilon)}\epsilon$ is
increasing in $\epsilon$ and converges to $\gamma_2=2^{\frac{\rho}
{\rho-1}}\gamma_1$  as $\epsilon\to 0$.

\smallskip
\subsubsection{Unbounded potentials with double exponential type tails}
Here we consider the second universality class so that
$H$ is regularly varying with index $\gamma=1$ and $k^*\in(0,\infty)$.
As shown in \cite{hkm} in proposition 1.1, this
is equivalent to the
 existence of a constant $\rho\in (0,\infty)$ such that

\begin{equation}
\label{ass-h}
\lim_{t\to\infty}\frac{H(yt)-yH(t)}{t}=\rho y\log y,
\end{equation}
for all $y\in (0,1)\cup(1,\infty)$ (this and $0<\rho<\infty$ is called
assumption {\bf (H)} in \cite{gm2}). Furthermore, under
assumption (\ref{ass-h}) it is true that,

$$
\lim_{t\to\infty}\frac{H(t)}{t}=\infty.
$$
\noindent This second
universality class  includes the case of unbounded
potentials which are double exponentially distributed with
parameter $\rho$, $0<\rho<\infty$,

$$
\mu[v(0)>x]=\exp\left\{-e^{x/\rho}\right\},
$$
for $x\in{\mathbb R}$.
We then have the following interesting
proposition.
\smallskip

\begin{proposition}
\label{double}
 Suppose that assumption (\ref{ass-h}) is satisfied for
some $\rho\in (0,\infty)$. Then assumption {\bf (RI)} is satisfied with
transition exponents

$$\gamma_1=\rho,\qquad\qquad
\gamma_2= 2\rho,$$
and growth exponent,

$$J(t)=t.$$
\end{proposition}

\smallskip

\noindent To prove proposition \ref{double}, we
quote  theorem 1.2 of \cite{gm2}, which shows that
under assumption (\ref{ass-h}) we have that,

$$
\langle m(0,t)\rangle=\exp\left\{H(t)-2d\kappa \chi\left(
\frac{\rho}{\kappa}\right)t+o(t)\right\},
$$
for $\kappa>0$, where $\chi(x):=\frac{1}{2}\inf_{p\in{\mathcal P}({\mathbb Z})}
[S(p)+\rho I(p)]$ for
$x\ge 0$, ${\mathcal P}({\mathbb Z})$ is the space of probability
measure on ${\mathbb Z}$, $S:{\mathcal P}({\mathbb Z})
\to [0,\infty)$ is the Donsker-Varadhan functional
defined by $S(p):=\sum_{x\in{\mathbb Z}}\left(
\sqrt{p(x+1)}-\sqrt{p(x)}\right)^2$, while $I:{\mathcal P}({\mathbb Z})
\to [0,\infty)$ is the entropy functional defined by
$I(p):=-\sum_{x\in{\mathbb Z}} p(x)\log p(x)$.
 On the other hand, from the previous
discussion we can  conclude that,
for every $\epsilon\ne 0$,

$$
 \frac{F_\epsilon(t)}{t}\sim\rho (1+\epsilon)
\frac{\log (1+\epsilon)}{\epsilon},
$$
and

$$
 \frac{F_\epsilon(2t)}{t}\sim 2\rho (1+\epsilon)
\frac{\log (1+\epsilon)}{\epsilon}.
$$
>From these limiting behaviors, we see that assumption
{\bf (MI)} is satisfied. Furthermore, form the fact that
$(1+\epsilon)\frac{\log (1+\epsilon)}{\epsilon}$ is increasing
for $\epsilon$ small enough and converges to $1$ as $\epsilon\to 0$,
we obtain the transition exponents at $\rho$ and $2\rho$
of proposition \ref{double}.
\smallskip

\subsubsection{Almost bounded potentials} We now focus on the
third universality class, where $\gamma=1$ and $k^*=0$. As shown
in theorem 1.4 of \cite{hkm}, in this case it is true that,

\begin{equation}
\label{u-three}
\log\langle m(0,t)\rangle \sim \frac{H(t\alpha_t^{-d})}{\alpha_t^{-d}},
\end{equation}
where $\alpha_t:[0,\infty)\to [0,\infty)$ is the
so called {\it scaling function}, which is implicitly defined for all
$t>0$ sufficiently large by the equation,

\begin{equation}
\label{imp}
\frac{k(t\alpha_t^{-d})}{t\alpha_t^{-d}}=\frac{1}{\alpha_t^2}.
\end{equation}
As shown in proposition 1.2 of \cite{hkm}, this function is unique
up to asymptotic equivalence. Furthermore, for the third universality
class, it is a slowly varying function (regularly varying of index $0$).

\smallskip

\begin{proposition}
\label{third}
 Suppose that $H$ is regularly varying of index $1$
and that (\ref{p11})
 is satisfied for $k(t)$ such that $k(t)\ll t$.
 Then assumption {\bf (RI)} is satisfied with
transition exponents

$$\gamma_1=\rho,\qquad\qquad
\gamma_2= 2\rho,$$
and growth exponent,

\begin{equation}
\label{jay}
J(t)=\frac{t}{\alpha_t^2}.
\end{equation}
\end{proposition}

\smallskip
\noindent The proof of proposition \ref{third} follows now
applying the definition of $\alpha_t$ through (\ref{imp}),
the asymptotics (\ref{u-three}) and (\ref{p11}).

A specific example of a distribution falling in the third universality
class is given by the {\it squared double exponential} law.
In other words, an effective potential $v(0)$ which is unbounded,
with law,

$$
\mu[v(0)>x]=\exp\left\{-e^{x^2}\right\},
$$
for $x\ge 0$. As it can be deduced from the discussion of example 1.4.3
of \cite{hkm}, in this case the scale function is given by $\alpha_t
\sim 2^{1/2}(\log t)^{1/4}$ and $k(t)\sim\frac{t}{2\sqrt{\log t}}$.
Furthermore, $\rho=1$ and the growth exponent can be chosen as
$J(t)\sim t/(2(\log t)^{1/2})$ (see table 2 of the introduction).

\smallskip

\subsubsection{Bounded potentials with Fr\'echet type
 tails} We now continue with an example falling in the
universality class $(4)$ of \cite{hkm}. We will consider
the case in which the essential supremum of $v(0)$ is $0$
and the tails are of Fr\'echet-type:
$\mu[v(0)>-x]=\exp\left\{-h(x^{-1}\right\}$ for $x>0$, with
$h\in R_\rho$ for some $0<\rho<\infty$. Using the
terminology of \cite{bbm} in the context of i.i.d. random exponentials,
this is {\it Case } {\bf A}.
To state appropriately this result, we will need to recall the work of Biskup and K\"onig \cite{bk}, who
studied the asymptotics of the annealed and quenched
first moments in
the case of product environments $\mu$ such that the cumulant
generating function $H(t)$ is in the so called $\gamma$-class
for some $\gamma\in [0,1)$. We say that {\it $H$ is in the
 $\gamma$-class} if the
essential supremum of $v(0)$ is $0$ and if there is a non-decreasing
function $\alpha_t\in (0,\infty)$ and a function $\tilde H:
[0,\infty)\to (-\infty,0]$, $\tilde H\ne 0$, such that

\begin{equation}
\label{scale}
\lim_{t\to\infty}\frac{\alpha_t^{d+2}}{t}H\left(\frac{t}{\alpha_t^d}y
\right)=\tilde H(t),
\end{equation}
for $y\ge 0$, uniformly on compact sets in $(0,\infty)$.
We will denote the function $\alpha_t$ the {\it scale function}.
It is not difficult to show, using the de Bruijn exponential
Tauberian theorem \cite{bgt}, that
if $v(0)$ is in case {\bf A} for some $\rho>0$, then $H$
is in the $\rho'$-class for $\rho'=\frac{\rho}{\rho+1}$, or
$\frac{1}{\rho'}-\frac{1}{\rho}=1$.
We can now sate the following proposition.

\smallskip
\begin{proposition}
\label{frechet}
 Suppose that the essential supremum
of the effective potential $v(0)$ is $0$ with
Fr\'echet-type tails $\mu[v(0)>-x]
=\exp\{-h(x^{-1})\}$ for $x>0$, and $h\in R_{\rho}$ for some
$0<\rho<\infty$. Then assumption {\bf (RI)} is satisfied with
transition exponents

$$\gamma_1=\left(\frac{1}{d+2+2\rho}\right)^2,\qquad\qquad
\gamma_2= 2^{1-\gamma_1}\gamma_1,$$
and growth exponent,

$$J(t)=\chi\frac{t}{\alpha^2_t},$$
for some constant $\chi\in (0,\infty)$.
\end{proposition}

\smallskip

\noindent As a consequence of proposition \ref{frechet} we
can now apply corollary \ref{corollaryr}, generalizing
case {\bf A} of theorems 2.1 and 2.2 of \cite{bbm} from
 $\kappa=0$ to  $\kappa\ge 0$. In contrast to
proposition \ref{weibull}, where there is no
change in the value of the transition exponents
$\gamma_1$ and $\gamma_2$ from $\kappa=0$ to $\kappa>0$, here
there is.

Let us now prove proposition \ref{frechet}.
In proposition 2.1 of \cite{bk}, it is shown
that whenever $\mu$ is in the $\rho'$-class,
then the scaling function $\alpha_t\in R_\nu$, for,

$$
\nu:=\frac{1-\rho'}{d+2-d\rho'}.
$$
Furthermore, theorem 1.2 of \cite{bk} states that when $\mu$
is in the $\rho'$-class,
 there
exists a $\chi\in(0,\infty)$ such that,

$$
\log\langle m(0,t)^\beta\rangle\sim -\chi \frac{\beta t}{\alpha_{\beta t}^2},
$$
for every $\beta\in(0,\infty)$.
Then condition {\bf (MI)} of theorem \ref{te1} is satisfied, and for
every $\epsilon\ne 0$, the  growth exponents  $F_\epsilon (t)$
and $F_\epsilon (2t)$  satisfy

$$
\frac{F_\epsilon(t)}{t/\alpha_t^2}
\sim\frac{(1+\epsilon)-(1+\epsilon)^{1-\nu^2}}\epsilon$$
and
$$
\frac{F_\epsilon(2t)}{t/\alpha_t^2}
\sim 2^{1-\nu^2}\frac{(1+\epsilon)-(1+\epsilon)^{1-\nu^2}}\epsilon.$$
As in the proof of
proposition  \ref{weibull}, we can show that
 that $\frac{(1+\epsilon)-(1+\epsilon)^{1-\nu^2}}\epsilon$ is
increasing in $\epsilon$ and converges to $\gamma_1=\nu^2$
as
$\epsilon\to 0$, and similarly
 $2^{1-\nu^2}\frac{(1+\epsilon)-(1+\epsilon)^{1-\nu^2}}\epsilon$ is
increasing in $\epsilon$ and converges to $\gamma_2=2^{1-\nu^2}
\gamma_1$  as $\epsilon\to 0$.

\smallskip

\subsection{The critical regime} For the examples discussed in the
four previous sections, it is possible to obtain the following
improvement of theorem \ref{te1}.

\begin{theorem}
\label{te2} Consider a reaction-diffusion process
with initial law $\delta_0$ and
field $w=(v_+,v_-)$ distributed according to a product measure
$\mu\in{\mathcal P}(W)$. Then the following statements are
satisfied.

\begin{itemize}

\item[(i)] {\bf Weibull-type}. Assume that $\mu$ has Weibull-type tails
so that
$\mu[v(0)>x]=\exp\left\{-h(x)\right\}$ for $x>0$, $h\in R_\rho$ for
$1<\rho<\infty$. Then if
$d\log L(t)\le \gamma H(t)$ eventually in $t$, and $0<\gamma<\gamma_1=
\frac{1}{\rho-1}$,
we have that for every $\delta>0$ in $\mu$-probability,

$$
\frac{m_L(0,t)}{e^{ (a_W(\gamma)+\delta) H(t)}}\ll 1,
$$
where

$$
a_W(\gamma):=\frac{\rho}{\rho-1}[(\rho-1)\gamma]^{1/\rho}-\gamma.
$$

\item[(ii)] {\bf Double exponential-type}. Assume that $\mu$
satisfies assumption (\ref{ass-h}) for some
constant $\rho\in (0,\infty)$. Then if,
$d\log L(t)\le \gamma t$ eventually in $t$,
and $\gamma<\gamma_1=\rho$, we have that for every $\delta>0$ in $\mu$-probability,

$$
\frac{m_L(0,t)}{
\exp\left\{ \frac{H((a_D(\gamma )+\delta)t)}{a_D(\gamma)+\delta}\right\} }\ll 1,
$$
where

$$
a_D(\gamma):=\gamma e^{\frac{1}{\rho}(\gamma-\rho)}.
$$

\item[(iii)] {\bf Almost bounded potentials}. Assume that  $\mu$
is such that
 $H$ is regularly varying of index $1$
and that (\ref{p11})
 is satisfied for $k(t)$ such that $k(t)\ll t$.
 Consider the growth exponent $J(t)$ as defined in
display (\ref{jay}). Then if,
$d\log L(t)\le \gamma J(t)$ eventually in $t$,
and $\gamma<\gamma_1=1$, we have that for every $\delta>0$
 in $\mu$-probability,

$$
\frac{m_L(0,t)}{
\exp\left\{ \frac{H((a_A(\gamma )+\delta)t)}{a_D(\gamma)+\delta}\right\} }\ll 1,
$$
where

$$
a_A(\gamma):=\gamma e^{\frac{1}{\rho}\left(\gamma-\rho\right)}.
$$

\item[(iv)] {\bf Fr\'echet-type}.
Assume that $\mu$ is such that  $essup v(0)=0$ and is of
Fr\'echet-type so that $\mu[v(0)>-x]=\exp\left\{-h(x^{-1})\right\}$
for $x>0$, and $h\in R_\rho$ for some $\rho\in (0,\infty)$.
Then if
$d\log L(t)\le \gamma J(t)$ eventually in $t$,
and $\gamma<\gamma_1=\nu^2$, we have that for every $\delta>0$
in $\mu$-probability,

$$
\frac{m_L(0,t)}{e^{ -(a_F(\gamma )-\delta)J(t)}}\ll 1,
$$
where

$$
a_F(\gamma):=(1-\nu^2)\left(\frac{\gamma}{\nu^2}\right)^{-
\frac{\nu^2}{1-\nu^2}}+\gamma.
$$

\end{itemize}

\end{theorem}

\noindent Let us remark that part $(iv)$ of theorem \ref{te2} includes
as a particular case, {\bf Case 3} of part $(i)$ of theorem 2 of \cite{bmr}.
We believe that the four functions $a_W, a_D, a_A, a_F$ are sharp,
in the sense that the quantities of the four parts of theorem
\ref{te2} diverge if the sign of $\delta$ is changed.
Also, these four functions have as maximum value $1$,
which is reached at $\gamma_1$.

In the special case in which $\mu$ is a double exponential law
so that $\log \mu[v(0)>x]=-e^{x^2}$
in part $(ii)$ of theorem \ref{te2}, the function
$\exp\left\{ \frac{H((a_D(\gamma )+\delta)t)}{a_D(\gamma)+\delta}\right\}$
takes the form
$\exp\left\{ t\log\left(\frac{(a_D(\gamma )+\delta)t}{e}\right)\right\}$,
showing how the transition mechanism takes place at a logarithmic
order in the exponent in contrast to the polynomial one
of parts $(i)$ and $(iv)$.
The whole picture suggested by
theorem \ref{te2} seems to indicate the presence of a phase transition
type behavior, as it is found in some mean field statistical mechanics
models like the Random Energy Model \cite{d}, \cite{bkl}.
 Indeed, when combined with
part $(i)$ of theorem \ref{te1}, we conclude that
$\log m_L(0,t,w)\sim \bar a(\gamma) J(t)$, where $\bar a(\gamma)$ equals
$a_W, a_D, a_A$ or $a_F$ depending on the potential for $\gamma<\gamma_1$,
while $\bar a(\gamma )=1$ for $\gamma >1$.
Thus, there is non-analyticity at $\gamma=1$ of a quantity playing the role of a ``free
energy''.

\medskip
\section{The  conditions for no explosion}
In this section we will prove proposition \ref{non-explosion}.
Since the  initial conditions $\nu\in{\mathcal P}(V)$
are concentrated on configurations with a finite number of particles,
and by translation invariance of the dynamics of the reaction-diffusion
process on ${\mathbb Z}^d$,
note that it is enough to consider the case where $\nu=\delta_0$.

\smallskip

\subsection{Preliminary lemmas}
Let us consider the reaction-diffusion process at scale $n$
with field $w$ satisfying condition  (\ref{assumet}),
 and with the initial condition $\delta_0$. Define the
quantities,

$$\zeta^n(t):=\sum_{x\in\Lambda_n}\eta^n(t,x),$$
representing the total number of particles produced at time $t$
and,

$$\bar\zeta^n(t):=\sum_{x\in\Lambda^c_n}(\eta^n(t,x)-\eta^n(0,x)),$$
 representing the total number of particles which have touched
$\Lambda_n^c$ in the time interval $[0,t]$. From this definition,
we can conclude that for $m>n\ge 1$ it is true that,

\begin{equation}
\label{true}
\zeta^m(t)=\zeta^n(t)+\sum_{k=1}^{\bar\zeta^n(t)}\zeta^m_{x_k}(t-\tau_k),
\end{equation}
where for $ 1\le k\le \bar\zeta^n(t)$,  $x_k\in\delta\Lambda_n$ is the
set of exit sites from $\Lambda_n$ of the random walks which have
touched the set $\Lambda_n^c$ in the time interval
$[0,t]$,
 $ 0\le\tau_k\le t$  the exit times
of each one of these random walks  and $\{\zeta^m_{x_k}(s):s\ge 0\}$ is a set
of independent processes such that  $\zeta^m_{x_k}(s)=0$ for $s<0$,
 $\zeta^m_{x_k}(s)=\delta_{x_k}$
 for $s=0$ while
$\{\zeta^m_{x_k}(s):s\ge 0\}$ has the law $P^{m,w}_{x_k}$.
Let us now define for each $n\in{\mathbb N}$ the maximum value
of the field $v_+$ on the box $\Lambda_n$ by $v_n:=\max_{x\in\Lambda_n}v_+(x)$.

\medskip

\begin{lemma}
\label{estimates}
 Consider the reaction-diffusion process  at scale
$n$ with field $w$ and initial condition $\delta_0$. Then,

\begin{itemize}

\item[(i)] For every $t\ge 0$ and $n\ge 1$,
$$E^{w}_0\left[ \zeta^n(t)\right]\le \exp\left\{v_n t\right\}.
$$

\item[(ii)] For every $t\ge 0$ and $n\ge 2\kappa t$,

$$E^{w}_0\left[ \bar\zeta^n(t)\right]\le 4d\exp\left\{\left(v_n-2\kappa\right) t
-n\log\left(\frac{n}{2e\kappa t}\right)\right\}.
$$

\item[(iii)] For every $t\ge 0$ and $m>n\ge 1$,

$$E^{w}_0\left[ \zeta^m(t)\right]\le
E^w_0\left[\zeta^n(t)\right]+E^w_0\left[\bar\zeta^n(t)\right]E^w_0\left[\zeta^{2n+m}(t)\right].
$$
\end{itemize}

\end{lemma}

\medskip

\noindent The following elementary lemma will be used to prove
lemma \ref{estimates}.
We will need to define,

\begin{equation}
\label{rate-function}
I(y):=\sup_{\lambda\in{\mathbb R}}\left\{\lambda y
-(\cosh\lambda -1)\right\}=y\sinh^{-1}y-\sqrt{1+y^2}+1,
\end{equation}
Note that $I:[0,\infty)\to [0,\infty)$ is one to one and that
$I(x)>0$ for $x>0$.

\smallskip

\begin{lemma}
\label{poisson}
 Let $\{X_t:t\ge 0\}$ be a simple symmetric continuous time random walk
on ${\mathbb Z}$ of total jump rate $2\kappa>0$ starting
from $0$. For a
non-negative real $x$ define $\tau_x:=\inf\{t\ge 0:|X_t|\ge x\}$
as the first exit time of this random walk from the interval $\Lambda_x$.
Then, if $P$ is its law, we have

\begin{equation}
\label{realx}
P[\tau_x<t]\le 4\exp\left\{-2\kappa t I\left(\frac{x}{2\kappa t}\right)
\right\}
\le
4\exp\left\{-2\kappa t-x\log \left(\frac{x}{2e\kappa t}\right)
\right\},
\end{equation}
the second inequality being satisfied only for $x\ge 2\kappa t$.

\end{lemma}
\smallskip

\noindent The proof of lemma \ref{poisson} is a simple large deviation
estimate and will be omitted.
\smallskip

\begin{proof}[Proof of lemma \ref{estimates}]{\it Part (i)}. Let us remark that for every bounded non-decreasing function $f:{\mathbb N}
\to {\mathbb R}$ which is eventually constant we have,

$$L_n f(\zeta^n)\le v_n \zeta^n (f(\zeta^n+1)-f(\zeta^n)).$$
For a natural $N\ge 1$ fixed, choose $f(m)=m\land N$. Using the fact that
$\zeta^n(0)=1$, we then conclude that,

\begin{equation}
\label{dethat}
\zeta^n(t)\land N-1-v_n\int_0^t \zeta^n(s) \theta_{[0,N-1]}(\zeta^n(s))ds,
\end{equation}
where for $A\subset {\mathbb R}$, $\theta_A$ is the indicator
function of the set $A$, is a super-martingale. Hence, since the integrand
of the integral in (\ref{dethat}) is a positive function, by Fubini's
theorem,

$$E^{w}_0\left[\zeta^n(t)\theta_{[0,N-1]}(\zeta^n(t))\right]\le
1+v_n\int_0^tE^{w}_0\left[\zeta^n(s)\theta_{[0,N-1]}(\zeta^n(s))\right]ds. $$
Therefore, by Gronwall's lemma,

$$E^{w}_0\left[\zeta^n(t)\theta_{[0,N-1]}(\zeta^n(t)\right]\le \exp\left\{v_n t\right\}.$$
Taking the limit when $N\to\infty$ and using the monotone convergence theorem
we conclude the proof of part $(i)$ of the lemma.

\medskip
\noindent {\it Part (ii)}. Let us note the following identity,

\begin{equation}
\label{identity}
E^w_0\left[\bar\zeta^n(t)\right]\le d e^{v_nt}P[\tau_n<t],
\end{equation}
where $\tau_n$ is the first exit time of a simple symmetric
continuous time random walk of total jump rate $2d\kappa$, starting
from the origin $0$,
from the box $\Lambda_n$ and $P$ is its law.
Now, from the second inequality of display (\ref{realx}) of
lemma \ref{poisson} applied to each of the $d$
coordinates of such a random walk, we conclude that,

$$P[\tau_n<t]\le 4d \exp\left\{-2\kappa t-n\log\left(\frac{n}{2e\kappa t}\right)\right\},$$
substituting the corresponding expression back
in (\ref{identity}) and using part $(i)$ of the proposition we conclude the proof of the lemma.
\medskip

\noindent {\it Part (iii)}. From (\ref{true})  we have that,

$$
E^w_0[\zeta^m(t)]=
E^w_0[\zeta^n(t)]+E^w_0\left[\sum_{k=1}^{\bar\zeta^n(t)}E^w_0[\zeta^m_{x_k}(t-\tau_k)|\bar\zeta^n(t)]\right].
$$
But,
$E^w_0[\zeta^m_{x_k}(t-\tau_k)|\bar\zeta^n(t)]\le E^w_0[\zeta^{2n+m}(t)]$,
which concludes the proof.
\end{proof}

\smallskip

\subsection{Proof of proposition \ref{non-explosion}}
 We will now prove proposition \ref{non-explosion} with the
help of lemma \ref{estimates}. Let $\delta=\frac{1}{5}$ and choose $N$ so that $v_n\le \frac{\delta}{t} \left(
n\log \left(\frac{n}{2e\kappa t}\right)-4n\log 4\right)$
 whenever $n\ge N$. By part $(i)$ of the lemma we have that,

$$E_x^w[\zeta^n(t)]\le
\exp\left\{\delta\left(n\log\left(\frac{n}{2e\kappa t}\right)-4n\log 4\right)\right\},$$
while by part $(ii)$ we have,

$$E_x^w[\bar\zeta^n(t)]\le 4d\exp\left\{-\left(1-\delta\right)
n\log \left(\frac{n}{2e\kappa t}\right)
-4\delta n\log 4\right\},$$
whenever $n\ge N$. Choosing $m=2n>N$ in part $(iii)$ of the same lemma, it follows that,
$$E_x^w[\zeta^{2n}(t)]\le A_n+B_n E_x^w[\zeta^{4n}(t)],$$
where $A_n:=e^{ \delta n\log \left(\frac{n}{2e\kappa t}\right)}$ and
$B_n:=4de^{-\left(1-\delta\right)n\log \left(\frac{n}{2e\kappa t}\right)-4\delta n\log 4}$.
Repeating the bound for $2n$ and $m=4n$ and substituting back we get,
\begin{equation}
\nonumber
E_x^w[\zeta^{2n}(t)]\le A_n+B_nA_{2n}+B_nB_{2n} E_x^w[\zeta^{8n}(t)].
\end{equation}
Now, by induction on $m$, we get that,

\begin{equation}
\label{get-that}
E_x^w[\zeta^{2n}(t)]\le\sum_{k=0}^{m-1}c_k+c_m\frac{E_x^w\left[\zeta^{n2^{m+1}}(t)\right]}{A_{n2^m}},
\end{equation}
where $c_0:=A_n$ and for $k\ge 0$, $c_{k+1}:=c_k\frac{B_{n2^k}A_{n2^{k+1}}}{A_{n2^k}}$. Now,

$$\frac{B_{n2^k}A_{n2^{k+1}}}{A_{n2^k}}\le 2d e^{-(1-2\delta)n2^k\log\left(\frac{n2^k}{2e\kappa t}\right)}.$$
Hence, by d'Alambert test and the fact that $1-2\delta>0$
 we know that the series $\sum_{k=0}^\infty c_k$ is convergent. On the other hand we have,

$$
c_m\frac{E_x^w\left[\zeta^{n2^{m+1}}(t)\right]}{A_{n2^m}}
\le c_{m-1}e^{\delta n 2^m \log\left(\frac{n 2^m}{2e\kappa t}\right)
-(1-2\delta)n2^{m-1}\log\left(\frac{n2^{m-1}}{2e\kappa t}\right)},$$
which tends to $0$ since $2\delta<1-2\delta$ and $c_m<1$ for $m$ large enough. Taking the limit
when $m\to\infty$, then when $n\to\infty$ and using the monotone convergence
theorem  in inequality (\ref{get-that}), we deduce that

$$
E_x^w[\zeta(t)]\le\sum_{k=0}^{\infty}c_k<\infty.
$$

\medskip

\section{Moment and correlation estimates}
\label{mce}
\smallskip

Here we will obtain some important bounds for the large time
asymptotic behavior of the field of quenched
 first moments $\{m(x,t,w):x\in{\mathbb Z}^d\}$,
the annealed first moment field $\{\langle m(x,t)\rangle :x\in{\mathbb Z}^d\}$,
and their correlations. In the first subsection, we will
prove proposition \ref{srel2}.

\smallskip

\subsection{Proof of proposition \ref{srel2}}
Our first lemma  states a  useful super-additivity and
convexity property of
the cumulant generating function of the random
variable $v(0)$.

\smallskip
\begin{lemma}
\label{variable}
Consider the cumulant generating function
$H(t):[0,\infty)\to {\mathbb R}$ of the random
variable $v(0)$, defined in display (\ref{gen-cum}). Then,
the following statements are true.

\begin{itemize}
\item[(i)] $H$ is super-additive. In other words, if $t_1,\ldots,t_n$
are non-negative reals then,

$$H(t_1+\cdots +t_n)\ge H(t_1)+\cdots +H(t_n).$$

\item[(ii)] Assume that
$\lim_{t\to\infty}tH''(t)=\infty$. Then,
for every $\alpha>1$,

\begin{equation}
\label{every}
\lim_{t\to\infty}\frac{H(\alpha t)-\alpha H(t)}{t}=\infty.
\end{equation}

\end{itemize}

\end{lemma}

\begin{proof}{\it Part (i)}. Let $t_1,t_2\ge 0$.
>From H\"older's inequality, we have that
$\phi(t_1)\phi(t_2)\le \phi(t_1+t_2)$. Hence, $H(t_1+t_2)\ge
H(t_1)+H(t_2)$. By induction on $n$ we conclude the proof.

\smallskip

\noindent {\it Part (ii)}. From the assumption, note that
we can write $H''(t)=\frac{f(t)}{t}$, where $\lim_{t\to\infty}f(t)=\infty$.
Integrating the function $H''$ from $t$ to $\alpha t$, it follows
that,

\begin{equation}
\label{hderivative}
H'(\alpha t)-H'(t)\ge \inf_{s\ge t}f(s)\cdot\log\alpha.
\end{equation}
Integrating again we obtain for $u>t$ that,

$$H(\alpha t)-\alpha H(t)\ge
(t-u)
{\rm inf}_{s\ge u}
f(s)
\cdot\alpha\log\alpha+c(u),$$
where $c(u)=H(\alpha u)-\alpha H(u)$. Dividing by $t$, taking the limit
when $t\to\infty$ and then the limit when $u\to\infty$, we obtain
(\ref{every}).

\end{proof}

\smallskip
\noindent Let us now prove parts $(i)$ and $(ii)$ of
proposition \ref{srel2}. Note that,

\begin{equation}
\nonumber
\frac{\partial G}{\partial t}=
\frac{(1+\theta) H'((1+\theta)t)
-(1+\theta)H'(t)}{\theta}.
 \end{equation}
Then, inequality (\ref{hderivative}) implies part $(i)$ of the lemma.
On the other hand we have that,

\begin{equation}
\nonumber
\frac{\partial G}{\partial \theta}=
\frac{\theta t H'((1+\theta)t)
-H((1+\theta)t)+H(t)}{\theta^2},
 \end{equation}
By the mean value theorem there exists a $\bar \theta$ such that
$0<\bar\theta <\theta$ and $H((1+\theta)t)-H(t)=\theta tH'((1+\bar\theta)t)$ and
hence $\theta t H'((1+\theta)t)-H((1+\theta)t)+H(t)=
\theta t \left[H'((1+\theta)t)-H'((1+\bar \theta)t)\right]$, which
by inequality (\ref{hderivative}) is positive if $t\ge t_1$, where $t_1$
is independent of $\theta$ and $\bar\theta$. This proves
that $G_\theta(t)$ is monotone in  $\theta$ for $t\ge t_1$.

\smallskip
\noindent We now show that condition (\ref{condition3}) ensures
{\bf (SI)}. We will without loss of generality assume
that $\kappa>0$.
 Let us define real valued functions
$f,g$ by $f(x):=H(\theta t+xt)-H(xt)$ and $g(x):=2x$ for real $x$.
By the generalized mean value theorem applied in the interval $[1,1+\theta]$,
there exists a $\theta_1\in (0,\theta)$, such that
$\frac{f(1+\theta)-f(1)}{g(1+\theta)-g(1)}=\frac{f'(1+\theta_1)}
{g'(1+\theta_1)}$. In other words, the expression
$G_{2\theta}(t)-G_\theta (t)=\frac{H((1+2\theta)t)-2H((1+\theta)t)+H(t)}{2\theta}$, equals,

$$
\frac{t}{2}\left(H'((1+\theta +\theta_1) t)-H'((1+\theta_1)t)\right)
=\frac{\theta t^2}{2}H''((1+ \theta_1+\theta_2)t),
$$
where in the last equality we have applied the mean value theorem
and $\theta_2\in (0,\theta)$. It therefore follows that there
is a function $\bar\theta :[0,\infty)\to (0,2\theta)$ such that,

$$\frac{G_{2\theta}(t)-G_\theta (t)}{t}=\frac{\theta t}{2}
H''((1+\bar\theta (t))t).$$
Our hypothesis $\lim_{t\to\infty} tH''(t)=\infty$, shows that
the expression above tends to $\infty$ as $t\to\infty$.
The proof of  the case in which $\theta<0$   is similar
and the details will be omitted.

\noindent We now continue in subsection \ref{truncated}, defining
the truncated quenched first moments, and then describing the
parabolic Anderson equation satisfied by the quenched first moments
and the corresponding Feynman-Kac representations.

\medskip

\subsection{ Truncated quenched first moments}
\label{truncated}
 In the sequel, given a real
function $f(x)$ defined on ${\mathbb Z}^d$, we will define the discrete
Laplacian by,

\begin{equation}
\label{laplacian}
\Delta f(x):=\sum_{e\in{\mathbb Z}^d:|e|=1} \left(
f(x+e)-f(x)\right).
\end{equation}
\noindent Let us  now for each finite set $U\subset{\mathbb Z}^d$
and environment $w\in W$,
  define the field $w^{U}:=(v_-^{U},v_+)$ with
$v_-^{U}(x)=v_-(x)$ for $x\in U$, while
$v_-^{U}(x)=\infty$ for $x\notin U$.
We now define for $x\in{\mathbb Z}^d$ and $t\ge 0$,

$$\tilde m_U(x,t,w):=m(x,t,w^U).$$
As it will be seen later, this expression
satisfies the parabolic Anderson
equation with Dirichlet boundary conditions.
We will denote this quantity the {\it truncated quenched first moment}
on $U$ at time $t$ for a reaction-diffusion process
starting from $x$. Also, we will
call the set $\{\tilde m_U(x,t,w):x\in{\mathbb Z}^d\}$, the
{\it field of truncated first moments} on $U$ at time $t$.
Now, in the particular case
in which $U=\Lambda(x,r)$ for some $r>0$, we will use the
notation $\tilde m_r(x,t,w)$ instead of $\tilde m_U(x,t,w)$. We will refer
to this quantity as the {\it truncated quenched first moment} at scale $r$
at time $t$ for a reaction-diffusion process starting from site $x$.
Furthermore, we will call the sets $\{\tilde m_r(x,t,w):x\in{\mathbb Z}^d\}$,
the {\it field of truncated quenched first moments} at scale $r$ at time
$t$.

\medskip

\subsection{The parabolic Anderson equations}
 Here we will recall the moment equations
satisfied by the field of quenched first moments $\{m(x,t,w)\}$
and by the corresponding truncated fields.
Following \cite{gm},
we have the  proposition.

\begin{proposition} Let $U\subset{\mathbb Z}^d$ be a finite set
and $w\in W$ an environment.
Consider the field of quenched first moments $\{m(x,t,w):x\in{\mathbb Z}^d\}$
on ${\mathbb Z}^d$ at time $t$ and
the field of truncated quenched first moments
$\{\tilde m_U(x,t,w):x\in{\mathbb Z}^d\}$
on $U$ at time $t$. Then the following statements
are true.

\begin{itemize}

\item[(i)] The field of
truncated
quenched first moments $\{\tilde m_U(x,t,w):x\in{\mathbb Z}^d\}$
 of the total number of particles
at time $t$ on $U$, satisfies the equation,

$$\frac{\partial \tilde m_U}{\partial t}=\kappa\Delta \tilde m_U+v (x) \tilde m_U,\qquad
{\rm for}\quad x\in U\cap{\mathcal G}(w)^c$$
$$\tilde m_U(x,0,w)=1,\qquad\qquad {\rm for}\quad x\in {\mathbb Z}^d.$$
$$\tilde m_U(x,t,w)=0,\qquad\qquad {\rm for}\quad x\notin U\cap{\mathcal G}(w)^c, t>0.$$

\item[(ii)] The field of
quenched first moments $\{m(x,t,w):x\in{\mathbb Z}^d\}$
 of the total number of particles
at time $t$ on ${\mathbb Z}^d$, satisfies the equation,

$$\frac{\partial  m}{\partial t}=\kappa\Delta  m+v (x)  m,\qquad
{\rm for}\quad x\in {\mathcal G}(w)^c$$
$$ m(x,0,w)=1,\qquad\qquad {\rm for}\quad x\in {\mathbb Z}^d.$$
$$ m(x,t,w)=0,\qquad\qquad {\rm for}\quad x\notin {\mathcal G}(w)^c, t>0.$$

\end{itemize}
\end{proposition}
\begin{proof}
 Consider the family of functions
 $\{u_z(x,t):=E^w_x[z^{\zeta(t)}]\}$, parametrized by
complex $z$ such that $0<|z|\le 1$. It is easy to see that,

$$\frac{\partial u_z}{\partial t}=\kappa\Delta u_z+v_+(x)u_z^2
-(v_+(x)+v_-(x))u_z+v_-(x),$$

$$u_z(x,0)=z,$$
for $x\in{\mathcal G}(w)^c$, while $u_z(x,t)=1$ for $t\ge 0$ and
$x\in{\mathcal G}(w)$.
Differentiating the above equation with respect to $z$ we
obtain part $(i)$. A similar proof can be
carried out for part $(ii)$.

\end{proof}

\smallskip

\subsection{Bounds on the quenched first moments}
 We will now obtain upper and lower bounds for the
annealed moments of the quenched first moments.
Let us first recall two elementary inequalities. For $n$ natural, let
$a_1,\ldots ,a_n$ be arbitrary real numbers. Then, for $r\ge 1$,
we have Jensen's inequality,

\begin{equation}
\label{jensen}
\left|\sum_{i=1}^n a_i\right|^r\le n^{r-1}\sum_{i=1}^n|a_i|^r,
\end{equation}
while for $0\le r\le 1$ we have,

\begin{equation}
\label{inequality}
\left|\sum_{i=1}^n a_i\right|^r\le
\sum_{i=1}^n |a_i|^r.
\end{equation}
We will also need to introduce for $L\ge 0$ the notation,

\begin{equation}
\label{notation}
M_L:=\max_{x\in\Lambda_L}|v(x)|.
\end{equation}
Let us recall the following lemma, contained in the statement of theorem 2.1 of \cite{gm}.

\smallskip

\begin{lemma}
\label{flemma}
 Consider a finite subset $U\subset{\mathbb Z}^d$ and $\mu\in {\mathcal P}(W)$.
 Assume that $\mu$ satisfies condition {\bf (E)}.
Then, $\mu$-a.s. for every $x\in{\mathbb Z}^d$ and $t\ge 0$,
 the quenched
first moment $m(x,t,w)$ on ${\mathbb Z}^d$
at time $t$ starting from  site $x$ admits the Feynman-Kac representation

\begin{equation}
\label{kac1}
m(x,t,w)=E_x\left[e^{\int_0^t v(X_s)ds}
1(\tau_{{\mathcal G}(w)}>t)
\right],
\end{equation}
and the truncated quenched first moment $\tilde m_U(x,t,w)$ on $U$
at time $t$ starting from $x$ also,

\begin{equation}
\label{kac2}
\tilde m_U(x,t,w)=E_x\left[ e^{\int_0^t v(X_s)ds}1(\tau_{U^c\cup{\mathcal G}(w)}>t)\right],
\end{equation}
where in both (\ref{kac1}) and (\ref{kac2}),
 $\{X_t:t\ge 0\}$ is a simple symmetric random walk of total
jump rate $2d\kappa$ starting from $x$, of law $P_x$, $E_x$
 is the expectation related to this law, and  for $A\subset{\mathbb Z}^d$
we define
$\tau_A:=\inf\{t\ge 0:X_t\notin A\}$.

\end{lemma}
\medskip

\noindent We can now apply lemma \ref{flemma} to obtain the first
estimates on the quenched first moments.

\smallskip

\begin{proposition}
\label{numbers}
 Consider a finite subset $U\subset{\mathbb Z}^d$ and $\mu\in {\mathcal P}(W)$.
 Assume that $\mu$ satisfies condition {\bf (E)}.
Then,

\begin{itemize}

\item[(i)] For each $x\in{\mathbb Z}^d$, $t\ge 0$ and $\beta >0$ there exists a constant
$C$ such that,

\begin{equation}
\label{sthat}
 e^{ H(\beta t)-2d\kappa t}\le
\left\langle m(x,t)^\beta
\right\rangle\le C (\kappa+t)^d e^{ H(\beta t)},
\end{equation}

\item[(ii)] For each $x\in U$, $t\ge 0$ and $\beta >0$,

\begin{equation}
\label{sthat2}
e^{ H(\beta t)-2d\kappa t}
\le
\left\langle \tilde m_U(x,t)^\beta
\right\rangle\le C (\kappa+t)^d e^{
H(\beta t)},
\end{equation}
where  $C$  is the constant of part $(i)$.

\item[(iii)]  For each $\beta >0$, $\gamma>0$ and $a>0$ we have that,

\begin{equation}
\label{havet}
\left\langle \left|m(x,t)-\tilde m_{\gamma(\kappa t)^a}(x,t)\right|^\beta\right\rangle\le
C(\gamma(\kappa t)^a+1)^d e^{-2 \beta\kappa t I\left(\frac{\gamma (\kappa t)^{a-1}}{2}\right)}
e^{H(\beta t)},
\end{equation}
for some constant $C>0$, where $I:[0,\infty)\to[0,\infty)$ is defined in
display (\ref{rate-function}).

\end{itemize}

\end{proposition}

\begin{proof} {\it Part (i).} The first inequality of display (\ref{sthat})
can be obtained from the Feynman-Kac representation (\ref{kac1}),
taking only into account the contribution of the path $X_s$ which stays
during the whole time interval $[0,t]$ at $x$ (page 637 of \cite{gm}).
 To prove the second inequality of (\ref{sthat}), let
us  note that by translation invariance it is
enough to prove the estimate for $\langle m(0,t)\rangle$.
On the other hand,

$$1(\tau_{{\mathcal G}(w)}>t)e^{\int_0^tv(X_s)ds}\le\sum_{n=0}^\infty
e^{M_{R_n} t}1(T_{n-1}\le t<T_n),$$
where $T_{-1}=0$, while for $n$ natural $T_n$ is the first exit time of the
random walk $\{X_t:t\ge 0\}$ from the box $\Lambda_{R_n}$,
with $R_n:= R_0 2^n$ and $R_0:=\max\{\kappa t,1\}$,
while $M_{R_n}:=\max_{x:\sn x\sn \le R_n} |v(x)|$ as defined
in display (\ref{notation}).
It follows that,

\begin{equation}
\label{fthat}
m(0,t,w)\le \sum_{n=0}^\infty e^{M_{R_n} t} P[T_{n-1}\le t].
\end{equation}
Let us also remark that since $\beta>0$, for each natural $n$ we have the
following inequality  which will be used soon,

\begin{equation}
\label{later}
\langle e^{\beta t M_{R_n}}\rangle\le (2(R_n+1))^d \exp\left\{H(\beta t)
\right\}.
\end{equation}
In fact, $\langle e^{\beta t M_{R_n}}\rangle \le\sum_{x\in [-R_n,R_n]^d}
\left\langle e^{\beta tM_{R_n}} 1(v(x)=M_{R_n})\right\rangle$.
Let us now consider the case $0<\beta\le 1$. Then, by an
 application of inequality (\ref{inequality}) to
estimate (\ref{fthat}), we conclude that,

\begin{equation}
\nonumber
m(0,t,w)^\beta\le e^{\beta t M_{R_0}}+\sum_{n=1}^\infty e^{\beta tM_{R_n}}
 P[T_{n-1}\le t]^\beta.
\end{equation}
Taking expectations on both sides of this inequality and applying
the estimate (\ref{later}) and
the second inequality of display (\ref{realx}) of
lemma \ref{poisson} for each of the $d$ coordinates of the
 underlying random walk, we obtain,

\begin{equation}
\nonumber
\left\langle m(0,t)^\beta
\right\rangle
\le 2^d(R_0+1)^d
e^{H(\beta t)}\left(1+
4d\sum_{n=1}^\infty 2^{nd} e^{-\beta R_0 2^{n-1}\log\left(
\frac{ R_0 2^{n-1}}{2e\kappa t}\right)}\right),
\end{equation}
 Now,
 since
$e^{-\beta R_0 2^{n-1}\log\left(
\frac{ R_0 2^{n-1}}{2e\kappa t}\right)}\le
e^{-\beta  2^{n-1}\log\left(
\frac{ 2^{n-1}}{2e}\right)}$,
and
$R_0+1\le 2R_0$,
 we see that there is a constant
$C$ such that inequality (\ref{sthat}) is satisfied for
$0<\beta\le 1$.

\smallskip
\noindent Let us now consider the case $\beta>1$. Let $\beta'>1$ be
defined by $\frac{1}{\beta}+\frac{1}{\beta'}=1$. Then, if we represent
the left hand side of (\ref{fthat}) as
$\sum_{k=0}^\infty e^{tM_{R_n}}P[T_{n-1}<t\le T_n]^{\frac{1}{\beta'}}
P[T_{n-1}<t\le T_n]^{\frac{1}{\beta}}$, by H\"older's inequality
we get that,

$$m(0,t,w)\le\left(\sum_{n=0}^\infty e^{t\beta M_{R_n}}
P[T_{n-1}<t\le T_n]\right)^{1/\beta}.$$
A computation similar to the case $0<\beta\le 1$ finishes the proof of
part $(i)$.

\smallskip

\noindent {\it Part (ii).} The first inequality
of display (\ref{sthat2}) can be deduced by an argument analogous
to the one leading to the first inequality of display (\ref{sthat}).
Now, note from the representations
(\ref{kac1}) and (\ref{kac2}) that
$\tilde m_U(x,t,w)\le m(x,t,w)$. Hence, the second inequality
of display (\ref{sthat2}) is a corollary of
the second inequality of display (\ref{sthat}).

\smallskip
\noindent{\it Part (iii).} Let us remark from the Feynman-Kac representations
(\ref{kac1}) for $m(0,t)$ and (\ref{kac2}) for $\tilde m_{\gamma
(\kappa t)^a}(0,t)$ that,

$$
m(0,t)-\tilde m_{\gamma (\kappa t)^a}(0,t)= E_x\left[e^{\int_0^t v(X_s)ds}1(\tau_{\Lambda_{\gamma (\kappa t)^a}^c}<t)
1(\tau_{{\mathcal G}(w)}>t)\right].
$$
Hence, as in part $(i)$ of the proof of this proposition,

\begin{equation}
\label{auxil}
|m(0,t)-\tilde m_{\gamma (\kappa t)^a}(0,t)|\le \sum_{n=1}^\infty e^{M_{R'_n}t}P[T'_{n-1}\le t< T'_n],
\end{equation}
where $T'_n$ is the first exit time of the random walk $\{X_t:t\ge 0\}$ from the box $\Lambda_{R'_n}$,
with $R'_n:=R'_0 2^n$, $R'_0:=\max\{\gamma (\kappa t)^a,1\}$ and
$M_{R'_n}:=\max_{x:||x||\le R'_n}|v(x)|$. Let us consider the case $0<\beta\le 1$. By
inequality (\ref{inequality}) we obtain that,

\begin{equation}
\nonumber
|m(0,t)-\tilde m_{\gamma (\kappa t)^a}(0,t)|^\beta\le \sum_{n=1}^\infty e^{t\beta M_{R'_n}}P[T'_{n-1}\le t]^\beta.
\end{equation}
 But,
as in (\ref{later}) we can conclude that
$\langle e^{\beta t M_{R'_n}}\rangle\le 2^d(R'_n+1)^d \exp\left\{H(\beta t)
\right\}$. Hence, from the first inequality
of display (\ref{realx}) of lemma \ref{poisson}
we see that,

\begin{eqnarray}
\nonumber
&\left\langle\left|m(0,t)-\tilde m_{\gamma (\kappa t)^a}(0,t)
\right|^\beta\right\rangle
\le d 4^\beta 2^d(R'_1+1)^de^{H(\beta t)}
 e^{-2 \beta \kappa t I\left(
\frac{R_0'}{2\kappa t}\right)}\\
&
\label{correc}
\times\sum_{n=1}^\infty
\frac{(2(R'_02^n+1))^d}{(2(R_0'2+1))^d}
e^{2\kappa t\beta\left(I\left(\frac{R_0'}{2\kappa t}\right)
-I\left(\frac{R_0'2^n}{2\kappa t}\right)\right)}.
\end{eqnarray}
Now, using the fact that $I$ is convex, non-decreasing and that
$I(0)=0$, we see that,

$$
I\left(\frac{R_0'}{2\kappa t}\right)
-I\left(\frac{R_0'2^n}{2\kappa t}\right)
\le -(2^n-1)I\left(\frac{\gamma (\kappa t)^{a-1}}{2}\right).
$$
Substituting this back in (\ref{correc}), we obtain
inequality (\ref{havet}) in the case $0<\beta\le 1$.
The case $\beta>1$ can be proved using similarly, using
H\"older's inequality as in part $(i)$.
\end{proof}
\smallskip
\noindent We end this section with important estimates involving the
growth functions.
Given a subset $U\subset{\mathbb Z}^d$, we define the
discrete Laplacian operator with effective potential
$v(0)=v_+(0)-v_-(0)\in[-\infty,\infty)$ by its action on functions
$f:{\mathbb Z}^d\to {\mathbb R}$, which vanish outside
$U_w:=U\cap{\mathcal G}^c(w)$ ($f(x)=0$ for $x\notin U
\cap{\mathcal G}^c(w)$) as,

$$(\Delta+v)f(x)=\sum_{e\in{\mathcal B}}(f(x+e)-f(x))+v(x)f(x),$$
where ${\mathcal B}$ is the set formed by the elements of the
basis of the free Abelian group ${\mathbb Z}^d$ and its inverses.
Defining $L^2(U_w):=\{f:\sum_{x\in{\mathbb Z}^d}f(x)^2<\infty,
f(x)=0\ {\rm if}\ x\notin U_w\}$,
we can check that $\Delta+v$ is a bounded self adjoint
operator on the Hilbert space $L^2(U_w)$ endowed with the
inner product $(f,g):=\sum_{x\in{\mathbb Z}^d}f(x)g(x)$.
We then define $\{\lambda_n(U,w):0\le n\le {\mathcal U}-1\}$ as the
set of eigenvalues of $\Delta+v$ in $L^2(U_w)$ in decreasing
order, where
${\mathcal U}$ is the total number of eigenvalues. We
will denote by $\psi^{U,w}_n$ the corresponding normalized
eigenfunctions. Let $r\ge 0$. In the  case
in which $U=\Lambda(x,r)$, we will employ the notation
$\{\lambda_n(x,r,w)\}$ instead of $\{\lambda_n(U,w)\}$
and $\psi^{x,r,w}_n$ instead of $\psi^{U,w}_n$.

 We can now state the following important lemma.

\begin{lemma}
\label{improvement}
 Consider the  quenched first moment $m(0,t,w)$.
Then, for every  $\beta> 0$, $a>1$
 and $t\ge 1$,
there is a constant $k_1(d,a,\beta)$ such that,

\begin{equation}
\label{parttwo}
\frac{k_1^{-1}}{ (\kappa t)^{da(\beta+1)}+1}
 \left\langle m(0,\beta t)\right\rangle \le
\left\langle
m(0,t)^\beta\right\rangle\le
k_1
((\kappa t)^{da(\beta+1)}+1) \langle m(0,\beta t)\rangle,
\end{equation}
and

\begin{equation}
\label{parttilde}
\frac{k_1^{-1}}{(\kappa t)^{da(\beta+1)}+1}
 \left\langle \tilde m_{(\kappa t)^a}(0,\beta t)\right\rangle \le
\left\langle
\tilde m_{(\kappa t)^a}(0,t)^\beta\right\rangle\le
k_1
((\kappa t)^{da(\beta +1)}+1) \langle \tilde m_{(\kappa t)^a}(0,\beta t)\rangle.
\end{equation}

\end{lemma}

\noindent \begin{proof} We will only prove (\ref{parttwo}),
being the proof of display (\ref{parttilde}) analogous.
Let us first show that
for every real $a>0$, $\beta>0$ and $t\ge 0$ there is a constant
$c(d,\beta,a)$ such that,

\begin{equation}
\label{spec-up}
m(x,t,w)\le ([2(\beta\kappa t)^a]+1)^{d/2} e^{t\lambda_0(x,(\beta \kappa t)^a,w)}
+4d e^{-2\kappa t
I\left(\frac{(\beta \kappa t)^{a-1}}{2}\right)}
e^{ tM_{(\beta\kappa t)^a}}.
\end{equation}
 Note that
 by inequality (\ref{auxil}) with $\gamma=\beta^a$,
we have that $m(x,t,w)\le 4d e^{-2\kappa t
I\left(\frac{(\beta\kappa t)^{a-1}}{2}\right)}e^{ tM_{(\beta\kappa t)^a}}+\tilde m_U(x,t,w)$,
with $U=\Lambda(x,(\beta \kappa t)^a,w)$.
We then need to estimate the
 truncated quenched first moment at scale $(\beta\kappa t)^a$,
$\tilde m_{(\beta\kappa t)^a}(x,t,w)$.
First remark the following expansion in terms of the
eigenvalues $\{\lambda_n(x,(\beta \kappa t)^a,w)\}$ and the corresponding eigenfunctions
$\{\psi^{x,(\beta\kappa t)^a,w}_n\}$,

\begin{equation}
\label{mc2}
\tilde m_{(\beta\kappa t)^a}(x,t,w)=\sum_{n=0}^{{\mathcal U}-1} e^{t\lambda_n(x,(\beta \kappa t)^a,w)}
\psi^{x,(\beta \kappa t)^a,w}_n(x)(\psi^{x,(\beta\kappa t)^a,w}_n,
 1(A)),
\end{equation}
where $A:=\Lambda(x,(\beta\kappa t)^a,w)$. Now, by the Cauchy-Schwartz inequality we see that the right-hand side
of equality (\ref{mc2}) is upper-bounded by $e^{t\lambda_0(x,(\beta\kappa t)^a,w)}$
$\times\left(\sum_{n=0}^{{\mathcal U}-1}
(\psi_n^{x,(\beta\kappa t)^a,w}, 1(\{x\}))^2
\sum_{n=0}^{{\mathcal U}-1}(\psi_n^{x,(\beta\kappa t)^a,w}, 1(A))^2
\right)^{\frac{1}{2}}$
which  in turn is upper-bounded by
$e^{t\lambda_0(x,(\beta\kappa t)^a,w)}\sqrt{|A|}$, where $ 1(\{x\})(y)$
equals $1$ if $y=x$ and $0$ otherwise. Using the fact that
 $|A|=|\Lambda (x,(\beta \kappa t)^a,w)|\le ([2(\beta\kappa t)^a]+1)^d$,
finishes the proof of (\ref{spec-up}).

\noindent Let us now show that for every finite subset $U\subset{\mathbb Z}^d$,
it is true that,

\begin{equation}
\label{spec-low}
\frac{1}{|U|}\sum_{z\in U}m(z,t)\ge\frac{1}{|U|}
e^{t\lambda_0 (U,w)}.
\end{equation}
First note the trivial inequality $m(z,t,w)\ge \tilde m_{U}(z,t,w)$.
We also have the expansion,
$\tilde m_{U}(z,t,w)=\sum_{n=0}^{{\mathcal U}-1}
 e^{\lambda_n(U,w)t}\psi_n^{U,w}(z)
\left(\psi_n^{U,w},
 1(U)\right)$.
Therefore we can  see
that,

$$\frac{1}{|U|}\sum_{z\in U}m(z,t)\ge\frac{1}{|U|}
e^{\lambda_0(U,w)t}(\psi_0^{U,w},
1(U))^2$$
$$\ge \frac{1}{|U|}e^{\lambda_0(U,w)t}
\sum_{z\in U}(\psi_0^{U,w})^2(z)=\frac{1}{|U|}
e^{\lambda_0(U,w)t},$$
where we have used in the second  inequality the fact that
$\psi^{U,w}_0(x)\ge 0$
and in the last inequality the normalization condition $\sum_{z\in U}
(\psi_0^{U,w})^2(z)=1$.

\noindent Let us now prove the second inequality
of (\ref{parttwo}). By Jensen's inequality (\ref{jensen}), in
the case $\beta\ge 1$, or inequality (\ref{inequality}), in
the case $0<\beta<1$,
applied to (\ref{spec-up}), note that for
some constant $c(a,d,\beta)$,

$$\left\langle m(x,t)^\beta\right\rangle\le
c ((\kappa t)^{da\beta/2}+1) \left\langle e^{\beta  t\lambda_0(x,(\beta \kappa t)^a,w)}\right\rangle
+4de^{-2\beta\kappa t
I\left(\frac{ (\beta\kappa t)^{a-1}}{2}\right)}
\left\langle e^{\beta t M_{(\beta\kappa t)^a}}\right\rangle.$$
Now, by the first inequality of display (\ref{sthat}) of  part $(i)$ of proposition \ref{numbers} and a computation similar to the one leading
to (\ref{later}), the second term of the
right hand side of the above inequality, is upper bounded by,

\begin{equation}
\label{eq-term}
4d 2^d((\beta\kappa t)^a+1)^d\exp\left\{-2\beta\kappa t
I\left(\frac{ (\beta\kappa t)^{a-1}}{2}\right)+2d\kappa t\right\}
\langle m(x, t)^\beta\rangle.
\end{equation}
On the other hand, by inequality (\ref{spec-low}) with
$U=\Lambda(0, (\kappa t)^a)$, we see that,
$\left\langle e^{\beta t\lambda_0(x,(\beta \kappa t)^a,w)}\right\rangle\le
 (2[(\beta \kappa t)^a]+1)^{d}\langle m(x,\beta t)\rangle $.
It follows that,

$$\left\langle m(x,t,w)^\beta\right\rangle\le
c ((\kappa t)^{da(\beta+1)}+1)\langle m(x,\beta t)\rangle,
$$
where we have used the fact that for $a>1$, the term (\ref{eq-term})
is negligible with respect to $\langle m^\beta\rangle$.
The second inequality of display (\ref{parttwo}) of part $(ii)$ now follows noting that for $a>1$ the second term in the
right-hand factor of the above inequality is negligible with respect to
the first term.

\noindent Let us now prove the first inequality of display (\ref{parttwo}).
In the case $0<\beta\le 1$, by inequality (\ref{inequality}) and the bound (\ref{spec-low}),
we have that,

\begin{equation}
\label{in}
\sum_{x\in U} m^\beta (x,t,w)\ge e^{\beta t \lambda_0(U,w)}.
\end{equation}
But when $\beta>1$, by Jensen's inequality  we have that
$\sum_{x\in U} m^\beta (x,t,w)\ge |U|^{-(\beta -1)}
\left(\sum_{x\in U} m (x,t,w)\right)^\beta$, and (\ref{in}) is still
satisfied.
Choosing $U=\Lambda (0,(\beta \kappa t)^a)$, and using translation invariance
we get,

$$\langle m(0,t)^\beta\rangle \ge (2[(\beta \kappa t)^a]+1)^{-d}
\left\langle e^{\beta t \lambda_0(0,(\beta \kappa t)^a,w)}\right\rangle.$$
Using again the bound (\ref{spec-up}), neglecting the second term,
we finish the proof in the case $0<\beta\le 1$.
\end{proof}
\smallskip

\noindent An important consequence of lemma \ref{improvement}, is
that it shows that assumption {\bf (MI)} implies the so
called {\it intermittency effect} \cite{gm}. Let us define for
$\theta\ne 0$, the functions,

$$
\bar F_\theta(t):=
    \frac{\log\langle m(0,t)^{1+\theta}\rangle-(1+\theta)
\log\langle m(0,t)\rangle}{\theta}.
$$
Note that by Jensen's inequality, $\bar F_\theta (t)\ge 0$.

\begin{corollary}
\label{cor-interm}
 Suppose that assumptions {\bf (E)} and
{\bf (MI)} are satisfied. Then, for every $\theta\ne 0$,

\begin{equation}
\label{interm}
\lim_{t\to\infty}\frac{\bar F_{2\theta}(t)-\bar F_\theta(t)}
{\theta \log(\kappa t+e)}
=\infty.
\end{equation}
\end{corollary}

\smallskip
\subsection{Correlation and variance estimates on the field of quenched first moments}
In order to prove part $(ii)$ of theorem \ref{te1}, it will be
important to have a control on the variance of the the
quenched first moments.
In the sequel of this paper, to avoid
heavy notation, we will use $m_a$ instead of $m_{(\kappa t)^a}$. Given $x,y\in{\mathbb Z}^d$ and $t\ge 0$, let us define,

$$c(x,y,t):=\left\langle m(x,t)m(y,t)\right\rangle-
\left\langle m(x,t)\right\rangle\left\langle m(y,t)\right\rangle,
$$
which we will call  the {\it correlation between sites $x$ and $y$ at time $t$}
of the field of quenched first moments. Similarly let us define for
$a>0$,

$$c_a(x,y,t):=\left\langle \tilde m_a(x,t)\tilde m_a(y,t)\right\rangle-
\left\langle \tilde m_a(x,t)\right\rangle\left\langle \tilde
m_a(y,t)\right\rangle,
$$
 the {\it correlation between sites $x$ and $y$ at time $t$} of the truncated
field of quenched first moments at scale $(\kappa t)^a$.
 Let us begin with the following lemma.

\smallskip

\begin{lemma}
\label{correlations} Let $t\ge 0$.
 Consider the fields of quenched first moments
$\{m(x,t,w):x\in{\mathbb Z}^d\}$
and truncated quenched first
moments at scale $(\kappa t)^a$, $\{\tilde m_a(x,t,w):x\in{\mathbb Z}^d\}$.
Then the following statements are true.

\begin{itemize}

\item[(i)] The sum of the correlations between site $0$ and the other sites
of the field of quenched first moments,
behaves asymptotically as $t\to\infty$ like
the sum of the correlations between site $0$ and the other sites
of the truncated field of quenched first moments at scale $(\kappa t)^a$.
In other words,

$$
\sum_{y\in{\mathbb Z}^d} c(0,y,t)
\sim
\sum_{y\in{\mathbb Z}^d} c_a(0,y,t).
$$

\item[(ii)] Let $\{U_t:t>0\}$ be a collection of subsets of the
lattice ${\mathbb Z}^d$ indexed by $t>0$. Assume that $|U_t|\sim |U_{t,(\kappa t)^a}|$
as $t\to\infty$, where $U_{t,r}:=\{x\in U_t: dist(x, U^c_t)\ge 2r\}$,
for $r>0$. Then,

\begin{equation}
\label{variance1}
Var_\mu\sum_{x\in\Lambda_{U_t}}m(x,t)\sim
|U_t|\sum_{y\in{\mathbb Z}^d} c(0,y,t),
\end{equation}

\begin{equation}
\label{variance2}
Var_\mu\sum_{x\in\Lambda_{U_t}}\tilde m_a(x,t)\sim
|U_t|\sum_{y\in{\mathbb Z}^d} c_a(0,y,t),
\end{equation}
and

\begin{equation}
\label{variance3}
Var_\mu\sum_{x\in\Lambda_{U_t}}m(x,t)\sim
Var_\mu\sum_{x\in\Lambda_{U_t}}\tilde m_a(x,t).
\end{equation}

\end{itemize}
\end{lemma}

\begin{proof}{\it Part (i)}. From the Feynman-Kac representation (\ref{kac1}) of
lemma \ref{flemma}, note that it is
possible to write,

$$m(x,t,w)=
E_x\left[e^{\sum_{z\in{\mathbb Z}^d} v(z) {\mathcal L}(z,t)}
1(A)
\right],
$$
where $A:={\{\tau_{{\mathcal G}^c(w)}>t\}}$,
 ${\mathcal L}(z,t):=\int_0^t \delta_z(X_s)ds$ is the local
time at the point $z$ of the random walk $\{X_t:t\ge 0\}$ starting
from $x$, and
$\delta_z:{\mathbb Z}^d\to \{0,1\}$ is the indicator function
of the set $\{z\}$.
>From here, using Fubini's theorem it follows that we have the following
representation for the correlations between site $x$ and $y$ of the
quenched field of first moments.

\begin{equation}
\label{tmoments}
c(x,y,t)=E_{x,y}\left[\left\langle e^{\sum_{z\in{\mathbb Z}^d}v
({\mathcal L}+\tilde{\mathcal L}) }
 1(A\cap \tilde A)\right\rangle
-
\left\langle e^{\sum_{z\in{\mathbb Z}^d}v
{\mathcal L} }
1(A)
\right\rangle
\left\langle e^{\sum_{z\in{\mathbb Z}^d}v
\tilde{\mathcal L} }
 1(\tilde A)
\right\rangle\right]
\end{equation}
where  $\tilde {\mathcal L}(z,t):=\int_0^t \delta_z(\tilde X_s)ds$ is the local
time at the point $z$ of the random walk $\{\tilde X_t:t\ge 0\}$, independent
of $\{X_t:t\ge 0\}$, starting
from $y$, with law $P_y$, and $\tilde A$ is an identical copy
of $A$, but defined in terms of the random walk $\{\tilde X_t:t\ge 0\}$.
 Furthermore, $E_{x,y}:=E_x\otimes E_y$, denotes
the expectation with respect to the law of the independent random walks
$\{X_t\}$ and $\{\tilde X_t\}$. Now, the expression (\ref{tmoments}) for
the correlations can be written in terms of the cumulant
generating function defined in display (\ref{gen-cum}),

\begin{equation}
\label{nas}
c(x,y,t)=E_{x,y}\left[ e^{\sum_{z\in{\mathbb Z}^d}
H({\mathcal L}+\tilde{\mathcal L}) }
-
 e^{\sum_{z\in{\mathbb Z}^d} H({\mathcal L}) }
 e^{\sum_{z\in{\mathbb Z}^d} H(\tilde{\mathcal L})
 }\right],
\end{equation}
where we have used the independence of the coordinates of the
effective potential $\{v(x):x\in{\mathbb Z}^d\}$ under $\mu$. Note that
the super-additivity of $H$ (part $(i)$ of lemma
\ref{variable}), implies that this expression is non-negative.
On the other hand, a reasoning similar to the one leading to the representation (\ref{nas}),
this time based on the Feynman-Kac representation (\ref{kac2}) of
lemma \ref{flemma},
enables us to deduce that,

\begin{equation}
\label{cethat}
\!\! c_a(x,y,t)\!\! =\!\! E_{x,y}\!\left[
 1(\tau_a>t) 1(\tilde\tau_a>t)\!\!
\left(\! e^{\sum_{z\in{\mathbb Z}^d}
H({\mathcal L}+\tilde{\mathcal L}) }\!\!
-
 e^{\sum_{z\in{\mathbb Z}^d} H({\mathcal L}) }
 e^{\sum_{z\in{\mathbb Z}^d} H(\tilde{\mathcal L})
 }\!\right)\!\right],
\end{equation}
where
$\tau_a:=\tau_{\Lambda(x,(\kappa t)^a)}=
\inf\{t\ge 0:X_t\not\in \Lambda(x,(\kappa t)^a)\}$
and
$\tilde\tau_a:=\tilde\tau_{\Lambda(x,
(\kappa t)^a)}:=\inf\{t\ge 0:\tilde X_t\not\in
\Lambda(y,(\kappa t)^a)\}$. From (\ref{nas})
and (\ref{cethat}) it follows that,

\begin{equation}
\label{wsthat}
\sum_{y\in{\mathbb Z}^d}c(0,y,t)\ge
\sum_{y\in{\mathbb Z}^d}c_a(0,y,t).
\end{equation}
But note that in reality, due to the independence between any pair of
truncated quenched first moments at time $t$
 at two points at a distance larger than $2(\kappa t)^a$, we have
$\sum_{y\in{\mathbb Z}^d}c_a(0,y,t)=\sum_{y\in\Lambda_{2
(\kappa t)^a}}c_a(0,y,t)$.
Furthermore,

\begin{equation}
\label{wehave}
\sum_{y\in{\mathbb Z}^d}c(0,y,t)=
\sum_{y\in\Lambda_{2(\kappa t)^a}}c(0,y,t)+
\sum_{y\not\in\Lambda_{2(\kappa t)^a}}c(0,y,t).
\end{equation}
Now, an application of the first inequality of display (\ref{sthat})
of part $(i)$ of proposition \ref{numbers} and of part $(ii)$ of the
same proposition, shows that since $a>0$, it is true that
 $\sum_{y\in\Lambda_{2(\kappa t)^a}}c(0,y,t)
 \sim \sum_{y\in\Lambda_{2(\kappa t)^a}}c_a(0,y,t)$.
And a second application of display (\ref{sthat}) and
 the first inequality in display (\ref{realx}) of
lemma \ref{poisson}, shows that $\sum_{y\not\in\Lambda_{2
(\kappa t)^a}}c(0,y,t)
\ll \sum_{y\in\Lambda_{2(\kappa t)^a}}c(0,y,t)$. This, together with inequality
(\ref{wsthat}), ends up the proof of part $(i)$ of lemma
\ref{correlations}.

\smallskip

\noindent {\it Part (ii)}. We will first prove (\ref{variance2}).
Remark that,

$$Var_\mu\sum_{x\in U_t} \tilde m_a(x,t)
=\sum_{x,y\in U_t} c_a(x,y,t)
=\sum_{x,y\in U_t:\sn x-y\sn\le 2(\kappa t)^a} c_a(x,y,t),$$
where we used as in the proof of part $(i)$, the
fact that $\tilde m_a(x,t,v)$ and $\tilde m_a(y,t,v)$
are independent for $\sn x-y\sn\ge 2(\kappa t)^a$. Now, the right-most member of
the above inequality is bounded by,
$|U_t|\sum_{y:\sn y\sn\le 2 t^a} c_a(0,y,t)$, which gives us the
inequality,

$$Var_\mu\sum_{x\in U_t}\tilde m_a(x,t)\le |U_t|\sum_{y\in{\mathbb Z}^d} c_a(0,y,t).$$
Similarly we can conclude that,

$$Var_\mu\sum_{x\in U_t}\tilde m_a(x,t)\ge |U_{t,(\kappa t)^a}|\sum_{y\in{\mathbb Z}^d} c_a(0,y,t).$$
This finishes the proof of the statement of display
(\ref{variance2}). The statement of
display (\ref{variance1}) now follows from display
(\ref{variance2}), part $(i)$ of
lemma \ref{correlations} proved above, and the inequalities,

$$Var_\mu\sum_{x\in U_t}\tilde m_a(x,t)\le
Var_\mu\sum_{x\in U_t} m_a(x,t)
\le |U_{t}|\sum_{y\in{\mathbb Z}^d} c(0,y,t).$$
Finally, note that (\ref{variance3}) is a trivial consequence of
(\ref{variance1}) and (\ref{variance2}).
\end{proof}

\medskip

\section{The Annealed and Gaussian regimes}

Here we will prove parts $(i)$ and $(ii)$ of theorem \ref{te1},
making use of the estimates obtained in the previous section,
and theorem \ref{te2}.
The main argument which will be used
to prove \ref{te1} is
a renormalization method, which we call  {\it partition
analysis}, as developed in \cite{bmr}.
In order
to present a self-contained proof, in the first subsection we
recall this technique. Then, in the second subsection, we
derive some important estimates via the partition analysis
technique, which will enable us to reduce the proofs to sums
of independent random variables. In the third subsection
 we prove the law of large numbers,
stated in display (\ref{int1}) of part $(i)$ of theorem \ref{te1}.
In the fourth subsection we will prove the negative part
(absence of a law of large numbers) stated in display (\ref{int2})
of part $(i)$ of theorem \ref{te1}. In subsection 5.5,
we will prove the central limit theorem stated in
display (\ref{int3}) of part $(ii)$ of theorem \ref{te1}.
Then, in subsection 5.6 we prove the absence of
a central limit behavior stated in display (\ref{int4})
of part $(ii)$ of the same theorem.
Finally, in subsection 5.7, we prove theorem \ref{te2}.

\medskip

\subsection{Partition Analysis} Here we shall follow closely
section 5.1 of Ben Arous, Molchanov and Ram\'\i rez \cite{bmr}.
 For a fixed natural $L$ we consider the box
 $\Lambda_L=\{x\in{\mathbb Z}^d: \sn x\sn\le L\}$.
 We will define two related but different kinds of partitions
 of $\Lambda_L$. The first one shows that
$\Lambda_L$ can be decomposed into disjoint {\sl partition boxes}
 $\{\Lambda'_{\bf i}:{\bf i}\in{\mathcal I}\}$, indexed by some set
${\mathcal I}$,
so that $\Lambda_L=\bigcup_{\bf i}\Lambda'_{\bf i}$.
The second one defines a partition of $\Lambda_L$ in a
{\sl strip set} and {\sl main boxes} $\{\Lambda''_{\mathbf i}:{\mathbf i}\in
{\mathcal I}\}$.  In the first case, the index set ${\mathcal I}$ will be
partitioned in disjoint subsets $\{{\mathcal I}_K:  K\in{\mathcal K}\}$,
where the cardinality of $\mathcal K$ is $2^d$,
in such a way that for each $K\in{\mathcal K}$ any pair of elements of the
collection of partition boxes
$\{\Lambda'_{\mathbf i}:{\mathbf i}\in{\mathcal I}_K\}$
is at a large Euclidean distance.  In the second partition case, it turns out that
the survival probabilities corresponding to sites in the strip set
have a total sum which is negligible, while the main boxes happen to be essentially
independent.
To proceed  we will need to introduce some notation defining the
corresponding scales and subsets.

Our first parameter is a natural number $L'$ smaller than or equal to $L$,
which will be called the {\sl mesoscopic
scale}.
By the division
algorithm, we know that there exist natural numbers $q$ and $\bar q$
such that $2L+1=qL'+\bar q$, with $0\le \bar q < q$. Note that this
last equation can be written in the form

\begin{equation}
\label{p1}
2L+1=\sum_{i=1}^qL'_i,
\end{equation} with
$L'_i=L'+\theta_{\bar q}(i)$ and $\theta_{\bar q}(i)=1$ for $i\le \bar q$ and
$\theta_{\bar q}(i)=0$ for $i>\bar q$. For our purposes, the relevant
fact is that $L'\le L'_i\le L'+1$. In the sequel, for any given
pair of real numbers $a,b$ we will use the notation $[a,b]_l$
for $[a,b]\cap{\bf Z}$. We now will subdivide the
box $[-L,L]_l$ in intervals according to equation (\ref{p1}).
Thus, we define $I_1:=\left[-L,-L+L'_1-1\right]_l$ and for $1<i\le q$ we
let $I_i:=\left[-L+\sum_{j=1}^{i-1}L'_i,-L+\sum_{j=1}^{i}L'_i-1\right]_l$.
Note that $I_q=\left[L-L'_q+1,L\right]_l$ and  $|I_i|=L'_i$.
 Next, we introduce a second parameter $r$ which is a natural number
smaller than or equal to $L'$. We will call $r$ the {\sl fine
scale}. Then, for
each $I_i$ we  define an interval $J_i$ such that
$J_i\subset I_i$, $|J_i|=L'-2r$ and the endpoints of $J_i$
are at a distance larger than $r$  to the endpoints of $I_i$.
To do so, first let $r_i:=r+\theta_{\bar q}(i)$. Then define
$J_1:=\left[-L+r,-L+L'_1-1-r_1\right]_l$ and for
$1<i\le q$ we let $J_i:=\left[-L+\sum_{j=1}^{i-1}L'_i+r,-L+\sum_{j=1}^iL'_i-1
-r_i\right]_l$.

We now proceed to define the partition in $\Lambda_L$ in partition boxes
and define the corresponding decomposition of the index set.
First we  define the set ${\mathcal I}:=\{1,2,\ldots ,q\}^d$, which
will correspond to the indexes parameterizing the sub-boxes.
 For a given element ${\bf i}\in{\mathcal I}$, of the form
${\bf i}=(i_1,\ldots ,i_d)$ with $1\le i_k\le q$, $1\le k\le d$, we define

\begin{equation}
\nonumber
\Lambda'_{\mathbf i}:=I_{i_1}\times I_{i_2}\times\cdots\times I_{i_d}.
\end{equation}
We will call such a set a {\sl partition box}.
 By definition
the cardinality $|\Lambda'_{\mathbf i}|$ of a partition
box satisfies,

\begin{equation}
\label{p2}
(L')^d\le |\Lambda'_{\mathbf i}|\le (L'+1)^d.
\end{equation}
Note also that the partition boxes defines a partition of $\Lambda_L$ so
that $\Lambda_L=\bigcup_{{\mathbf i}\in{\mathcal I}}\Lambda'_{\mathbf i}$ where
the union is disjoint.

Next we define a partition of the index set ${\mathcal I}$.
Consider the collection ${\mathcal K}$ of
subsets of $\{1,2,\ldots, d\}$. Note that $|{\mathcal K}|=2^d$.
 Now given $K\in{\mathcal K}$ we
 define ${\mathcal I}_K$ as the subset  of
${\mathcal I}$ having coordinates which are even for
$k\in K$ and odd for $k\notin K$. In other words, if we
define ${\mathbb E}$ as the set of even natural numbers
and ${\mathbb O}$ as the set of odd natural numbers then,

$${\mathcal I}_K:=\{{\mathbf i}=(i_1,\ldots ,i_d)\in{\mathcal I}:
i_k\in{\mathbb E}\ {\rm if}\ k\in K,\
i_k\in{\mathbb O}\ {\rm if}\ k\notin K,\ 1\le k\le d\}.$$
Note that $\{{\mathcal I}_K:K\in{\mathcal K}\}$ defines a
partition of ${\mathcal I}$ so that ${\mathcal I}=\bigcup_{K\in{\mathcal K}}
{\mathcal I}_K$ is a disjoint union.  Hence,

\begin{equation}
\label{disjoint-union}
\sum_{K\in{\mathcal K}}\sum_{{\mathbf i}\in{\mathcal I}_K}
\sum_{x\in\Lambda'_{\mathbf i}}f(x) =\sum_{x\in\Lambda_L}f(x),
\end{equation}
for every function $f:\Lambda_L\to{\mathbb R}$,
which is a consequence of the fact that
$\Lambda_L=\bigcup_{K\in{\mathcal K}}\bigcup_{{\mathbf i}\in{\mathcal I}_K}
\Lambda'_{\mathbf i}$ is a disjoint union.
 We will refer
to such a decomposition as the {\sl parity partition at scale $L'$} of
$\Lambda_L$. Furthermore,
given  $K\in{\mathcal K}$ and
any pair of boxes $\Lambda'_{\mathbf i}$ and $\Lambda'_{\mathbf j}$
with ${\mathbf i},{\mathbf j}\in{\mathcal I}_K$ and
${\mathbf i}\ne{\mathbf j}$ we have that,

\begin{equation}
\label{p3}
dist\left(\Lambda'_{\mathbf i},\Lambda'_{\mathbf j}\right)\ge L'.
\end{equation}
Here for any pair of subsets $A,B\subset{\mathbb Z}^d$
we define $dist(A,B):=\inf_{x\in A,y\in B}|x-y|$. In other
words (\ref{p3}) expresses the fact that the distance between any
pair of partition boxes with different sub-indexes in
${\mathcal I}_K$ is larger than or equal to $L'$. This
completes the description of the partition of $\Lambda_L$ in
partition boxes.

Next, we describe the partition of $\Lambda_L$ into the strip
set and main boxes. Given an ${\mathbf i}\in{\mathcal I}$
we let,

$$\Lambda''_{\mathbf i}:=J_{i_1}\times J_{i_2}\times\cdots\times J_{i_d}.$$
Such a box will be called a {\sl main box}. Its cardinality is
$|\Lambda''_{\mathbf i}|=(L'-2r)^d$. Now let,

$$S_L:=\Lambda_L-\bigcup_{{\mathbf i}\in\mathcal I}\Lambda''_{\mathbf i}.$$
Such a set will be called the {\sl strip set}. Note that $S_L$
and $\{\Lambda''_L:{\mathbf i}\in{\mathcal I}\}$ define a partition
of $\Lambda_L$. We will refer to such a partition as the
{\sl strip-box partition at scale $L'$} of $\Lambda_L$. We furthermore remark the
following cardinality estimate for the strip set which will
be useful later,

\begin{equation}
\label{p4}
\frac{|S_L|}{(2L+1)^d}\le \frac{\left((L'+1)^d-(L'-2r)^d\right)}{(L')^d},
\end{equation}
where we have used the fact that $|{\mathcal I}|=q^d$.

\smallskip

\subsection{Moment and decoupling inequalities}
Following \cite{bmr}, we recall here some standard inequalities
and derive inequalities involving sums of the quenched first moments
that will be necessary to perform the partition analysis.

Let us first recall the well-known inequality of von Bahr and Esseen
(page 82, exercise 2.6.20, of Petrov \cite{petrov}).

\smallskip

\begin{lemma} (von Bahr-Esseen). Let $X_1,\ldots, X_n$ be mean zero
independent random variables and $S_n:=\sum_{k=1}^n X_i$. Then if
$E$ denotes the expectation with respect to the joint law of these
random variables, and $1\le r\le 2$, it is true that

\begin{equation}
\label{bahr-esseen}
E|S_n|\le 2\sum_{k=1}^n E|X_k|^r.
\end{equation}

\end{lemma}

\smallskip

\noindent We continue with the following lemma (analogous to lemma 6 of \cite{bmr}),
which is  a consequence of (\ref{bahr-esseen}) and (\ref{inequality}).

\smallskip

\begin{lemma}
\label{ofand}
 Let $a>1$. Consider the field of truncated quenched first moments
$\{\tilde m_a(x,t,w):x\in{\mathbb Z}^d\}$ at scale $(\kappa t)^a$. Let $L(t), L'(t):
[0,\infty)\to {\mathbb N}$ be functions such that $(\kappa t)^a\le L'(t)\le L(t)$.
Then if $1\le r\le 2$, it is true that,

\begin{equation}
\label{truet}
\left\langle \left|\sum_{x\in\Lambda_L}(\tilde m_a-\langle \tilde m_a\rangle)
\right|^r\right\rangle\le 2(2L'+2)^{d(r-1)}(2L+1)^d
\left\langle \left|\tilde m_a-\langle \tilde m_a\rangle\right|^r\right\rangle.
\end{equation}

\end{lemma}

\begin{proof} The proof of this result is analogous to
that of lemma 6 of \cite{bmr}, requiring the use of the
decomposition (\ref{disjoint-union}) with $f(x)=
\tilde m_1 (x,t,w)-\langle\tilde m_1\rangle (x,t)$, von Bahr-Esseen
inequality (\ref{bahr-esseen}) and Jensen's inequality (\ref{jensen}).
\end{proof}

\medskip

\noindent Next, we have the following estimate analogous to
 lemma 7 of \cite{bmr}.

\medskip

\begin{lemma}
\label{followinge}
Let $a>1$.
Consider the field of quenched first moments
$\{m(x,t,w):x\in{\mathbb Z}^d\}$ and of truncated quenched
first moments $\{m_a(x,t,w):x\in{\mathbb Z}^d\}$ at scale $(\kappa t)^a$.
Let $L(t):[0,\infty)\to {\mathbb N}$ be such that $(\kappa t)^a\ll L$.
Then the following statements are true.

\begin{itemize}

\item [(i)] Assume  that,

$$
\left\langle\left|\frac{\sum_{x\in\Lambda_{L}} \tilde m_a}{(2L+1)^d
\langle \tilde m_a\rangle}-1\right|\right\rangle\ll 1.
$$
Then,

$$
\left\langle\left|\frac{m^L}{\langle m\rangle }- 1\right|\right\rangle\ll 1.
$$

\item[(ii)] Asymptotically as $t\to\infty$ we have,

\begin{equation}
\label{item}
\left\langle\left|\frac{\sum_{x\in\Lambda_L}(m-\langle m\rangle)}
{\sqrt{Var_\mu \sum_{x\in\Lambda_L} m}}-
\frac{\sum_{x\in\Lambda_L}(\tilde m_a-\langle \tilde m_a\rangle)}
{\sqrt{Var_\mu \sum_{x\in\Lambda_L} \tilde m_a}}\right|\right\rangle\ll 1,
\end{equation}
and

\begin{equation}
\label{item2}
\left\langle\left|\frac{\sum_{x\in\Lambda_L}(m-\langle m\rangle)}
{\sqrt{Var_\mu \sum_{x\in\Lambda_L} m}}-
\frac{\sum_{x\in\Lambda_L}(\tilde m_a-\langle \tilde m_a\rangle)}
{\sqrt{Var_\mu \sum_{x\in\Lambda_L}  m}}\right|\right\rangle\ll 1.
\end{equation}

\end{itemize}
\end{lemma}

\begin{proof} The case $\kappa=0$ is trivial, so we will assume
that $\kappa>0$.

\smallskip
\noindent {\it Part (i)}. A direct calculation shows us that,

$$
\left|\frac{m^L}{\langle m\rangle }
-\frac{\sum_{x\in\Lambda_{L}} \tilde m_a}{(2L+1)^d
\langle \tilde m_a\rangle}\right|
\le
\frac{\sum_{x\in\Lambda_L}\left|
m-\tilde m_a\right|}{(2L+1)^d\langle m\rangle}
+\frac{\langle \left|m-\tilde m_a
\right|\rangle}{\langle m\rangle}
\left|\frac{\sum_{x\in\Lambda_{L}}\tilde m_a}{(2L+1)^d
\langle\tilde m_a\rangle}\right|.
$$
Now, an application of inequality (\ref{havet}) of part $(ii)$
of proposition  \ref{numbers} and of the first inequality in
 display (\ref{sthat}) of part $(i)$ shows us that,
$\frac{\sum\langle|m-\tilde m_a|\rangle}{(2L+1)^d\langle m\rangle}
\le\varepsilon_1 (t)$, where
$\varepsilon_1(t):= d 4 2^d ((\kappa t)^a+1)^d
e^{-2\kappa t\left( I\left(\frac{(\kappa t)^{a-1}}{2}
\right)-d\right)}$.
 Similarly we have that $\frac{
\langle|m-\tilde m_a|\rangle}{\langle m\rangle}
\le\varepsilon_1 (t)$. By the triangle inequality it follows that,

\begin{equation}
\label{everyt}
\left\langle\left|\frac{m^L}{\langle m\rangle }
-\frac{\sum_{x\in\Lambda_{L}} \tilde m_a}{(2L+1)^d
\langle \tilde m_a\rangle}\right|\right\rangle\\
\le
2\varepsilon_1(t)
+\left\langle\left|\frac{\sum_{x\in\Lambda_{L}}\tilde m_a}{(2L+1)^d
\langle\tilde m_a\rangle}-1\right|\right\rangle.
\end{equation}
Now, for $a>1$, we have $I\left(\frac{(\kappa t)^{a-1}}{2}\right)\to\infty$.
Hence,
$\varepsilon_1(t)\to 0$ as $t\to\infty$. This clearly implies the
statement of part $(i)$ of the lemma.

\smallskip

\noindent {\it Part (ii)}. Let us first note that the left-hand side
of display (\ref{item}) is upper-bounded by,

\begin{equation}
\label{dby}
2\frac{\sum\limits_{x\in\Lambda_L}\left\langle\left|m-\tilde m_a\right|\right\rangle}
{\sqrt{Var_\mu\sum\limits_{x\in\Lambda_L}m}}
+\left\langle\frac{\sum\limits_{x\in\Lambda_L}|\tilde m_a-\langle\tilde m_a\rangle|}
{\sqrt{Var_\mu\sum\limits_{x\in\Lambda_L}\tilde m_a}}\right\rangle
\left(\frac{\sqrt{Var_\mu\sum\limits_{x\in\Lambda_L}\tilde m_a}}
{\sqrt{Var_\mu\sum\limits_{x\in\Lambda_L}
m}}-1\right).
\end{equation}
Now, by Cauchy-Schwartz inequality we have that
$\left\langle\frac{\sum_{x\in\Lambda_L}|\tilde
 m_a-\langle\tilde m_a\rangle|}
{\sqrt{Var_\mu\sum_{x\in\Lambda_L}\tilde m_a}}\right\rangle\le 1$.
On the other hand, by the assumption $ t^a\ll L$
we know that the hypothesis of part $(ii)$
of lemma \ref{correlations} is satisfied for $U_t=\Lambda_{L(t)}$ so that
the asymptotics (\ref{variance3}) of this lemma holds, and hence
the second term of (\ref{dby}) tends  to $0$ as $t\to\infty$.
Furthermore, Cauchy-Schwartz inequality and
part $(ii)$ of proposition \ref{numbers}, imply that
$\sum_{x\in\Lambda}\langle|m-\tilde m_a|\rangle\le
(2L+1)^d\varepsilon_2 (t) \sqrt{\langle m(0,t)^2\rangle}$, where
$\varepsilon_2(t):=d4 2^d((\kappa t)^a+1)^{d/2}
\exp\left\{-2\kappa tI\left(\frac{
(\kappa t)^{a-1}}
{2}\right)+2d\kappa t\right\}$. In addition, by part  $(ii)$ of lemma
\ref{correlations} we have
$Var_\mu\sum_{x\in\Lambda}m\sim (2L+1)^d\sum_{x\in{\mathbb Z}^d}c(0,y,t)$.
Now, $\sum_{x\in{\mathbb Z}^d}c(0,y,t)\ge \sqrt{\langle m^2\rangle
-\langle m\rangle^2}$. Thus,  the first term
of (\ref{dby}) is upper-bounded by a quantity which asymptotically
behaves as $t\to\infty$ like,

$$
\frac{\varepsilon_2(t)}{\sqrt{1-\langle m\rangle^2/\langle m^2\rangle}}.
$$
But by corollary \ref{cor-interm},
and the fact that $\bar F_1(t)\ge 0$, we conclude that
$\langle m\rangle^2/\langle m^2\rangle\ll 1$.
In brief,  the first term
of (\ref{dby}) is upper-bounded by a quantity which asymptotically
as $t\to\infty$ behaves like,
$\varepsilon_2(t)$. Obviously, when $a>1$ we have $\varepsilon_2(t)\ll 1$.
\end{proof}

\medskip

\subsection{The Annealed asymptotics}
Let us now prove the
law of large numbers stated in display (\ref{int1}) of
part $(i)$ of theorem \ref{te1}. To simplify the writing of
the expressions in the calculations, we will redefine $\epsilon$
as $2\epsilon$, assuming that

\begin{equation}
\label{ass1}
L(t)\ge \exp\left\{\frac{1}{d}F_{2\epsilon}
(t)\right\},
\end{equation}
for some $\epsilon>0$,
 and prove that then in $\mu$-probability it is true that,
\begin{equation}
\label{lln1}
\frac{ m^L(0,t,w)}{\langle m(0,t)\rangle}\sim 1.
\end{equation}
To do so we first remark
that by inequality (\ref{everyt}) of lemma \ref{followinge},
 it is enough to show that for  $a=3/2$,

\begin{equation}
\label{lln1.6}
\left\langle
\left|\frac{\sum_{x\in\Lambda_L}\tilde m_a}
{(2L+1)^d\langle \tilde m_a\rangle}-
1\right|^{1+\epsilon}\right\rangle\ll 1.
\end{equation}
Remark that the
right hand side of display (\ref{lln1.6}) can be rewritten as,

\begin{equation}
\label{lln1.7}
\left\langle\left|\frac{\sum_{x\in\Lambda_L}\tilde m_a}
{(2L+1)^d\langle \tilde m_a\rangle}-1\right|^{1+\epsilon}\right\rangle
=
\frac{\left\langle\left|\sum_{x\in\Lambda_L}\left( \tilde m_a-
\langle \tilde m_a \rangle\right)\right|^{1+\epsilon}\right\rangle}
{(2L+1)^{d(1+\epsilon)}\langle \tilde m_a\rangle^{1+\epsilon}}.
\end{equation}
At this point we  make use of the  parity partition
 decomposition for $\Lambda_L$ previously defined
to deal with the numerator of the right-hand side of
display (\ref{lln1.7}) via inequality (\ref{truet}) of lemma \ref{ofand}
 with $r=1+\epsilon$.
 We will chose a time dependent mesoscopic scale
$L'(t)=(\kappa t)^b$, where
$0<a<b$.
Therefore, by lemma \ref{ofand},
the right-hand side of equality (\ref{lln1.7}) is
upper-bounded by

\begin{equation}
\label{lln9}
\frac{ 2(2L'
+2)^{d\epsilon}\left\langle\left|
\tilde m_a-\langle\tilde m_a\rangle\right|^{1+\epsilon}\right\rangle}
{(2L+1)^{d\epsilon}\langle\tilde m_a\rangle^{1+\epsilon}}.
\end{equation}
 Now,
since for any non-negative reals $x,y$ we have $|x-y|^{1+\epsilon}
\le |x|^{1+\epsilon}+|y|^{1+\epsilon}$, it follows that
$\left\langle\left|\tilde m_a-\langle \tilde m_a\rangle
\right|^{1+\epsilon}\right\rangle\le
\left\langle \tilde m_a^{1+\epsilon}\right\rangle+\langle \tilde m_a
\rangle^{1+\epsilon}\le 2
\left\langle \tilde m_a^{1+\epsilon}\right\rangle$, where in the
last inequality we have used Jensen's inequality. Hence, since
$\tilde m_a\le m$,
the expression  of display (\ref{lln9}) is upper bounded by,

$$
\frac{ 2 (2L'+2)^{d\epsilon}
\left\langle m^{1+\epsilon}\right\rangle}
{(2L+1)^{d\epsilon}\langle\tilde m_a\rangle^{1+\epsilon}}.
$$
But by parts $(ii)$ and $(iii)$ of proposition \ref{numbers}, we
can replace in the denominator of the above
expression the term $\langle\tilde m_a\rangle$ by $\langle m\rangle$.
Thus, (\ref{lln9}) is upper-bounded by an expression which
is asymptotically equivalent to,

\begin{eqnarray}
\nonumber
&\ &\frac{2(2L'+2)^{d\epsilon+1}k_1 ((\kappa t)^{da(2+\epsilon)}+1)
 e^{
H_1((1+\epsilon)t)-
(1+\epsilon) H_1(t)}}
{(2L+1)^{d\epsilon}}\\
\label{boundedb}
&\ &\le
c((\kappa t)^{da(2+\epsilon)})+1)(L')^{d\epsilon}
e^{-\epsilon (F_{2\epsilon} (t)- F_{\epsilon}(t))},
\end{eqnarray}
where we have used the second inequality of display
(\ref{parttwo}) of lemma \ref{improvement} in the first inequality,
 and  assumption (\ref{ass1})
in the last inequality.
Note that when $\kappa =0$ the last expression reduces
to $e^{-\epsilon (G_{2\epsilon} (t)- G_{\epsilon}(t))}$.
 Then,  assumption {\bf (MI)}
shows that the first and the second terms
of the right-hand side of (\ref{boundedb}) tend to
$0$ as $t\to\infty$.

\medskip

\subsection{The non-Annealed asymptotics} In this subsection we will
prove the asymptotic behavior of display (\ref{int2}). Again, we will
redefine $\epsilon$ by $2\epsilon$, assuming that,

\begin{equation}
\label{nthat}
L(t)\le \exp\left\{ \frac{1}{d} F_{-2\epsilon}(t)\right\},
\end{equation}
 for some $\epsilon >0$.
Note that it will be enough to show that,

\begin{equation}
\label{wthat}
\left\langle\left|\frac{m_L(0,t)}{\left\langle m(0,t)\right\rangle}
\right|^{1-\epsilon}\right\rangle\ll 1.
\end{equation}
To do so, note from inequality (\ref{inequality}) that the left-hand
side of display (\ref{wthat}) is upper bounded by,

\begin{equation}
\label{bby}
(2L+1)^{d\epsilon}\frac{\left\langle m^{1-\epsilon}(0,t)\right\rangle}
{\left\langle m(0,t)\right\rangle^{1-\epsilon}}.
\end{equation}
Now, by the second inequality of display (\ref{parttwo}) of
lemma \ref{improvement} applied with $\beta=1-\epsilon$ to the numerator of
(\ref{bby}), we see that the left-hand side of display
(\ref{wthat}) is upper bounded by,

$$(2L+1)^{d\epsilon} k_1 ((\kappa t)^{da(2-\epsilon)}+1) e^{H_1((1-\epsilon)t)-(1-\epsilon)H_1(t)
}.$$
Finally, by assumption (\ref{nthat}), this expression is upper bounded
by,

$$c ((\kappa t)^{da(2-\epsilon)}+1) e^{-\epsilon (
F_{-\epsilon}(t)- F_{-2\epsilon}(t))}.$$
Now, the assumption {\bf (MI)}, implies that this
expression converges to $0$ as $t\to\infty$.

\medskip

\subsection{The Gaussian asymptotics} Here we prove the central limit
theorem stated in display (\ref{int3}) of part $(ii)$ of theorem \ref{te1}.
We will perform this time a strip-box partition of the box $\Lambda_L$
into the strip set $S_L$ and the main boxes. Let us fix $a>1$.
 Let us fix $b$ and
$c$ such that $c>b>a$, and
 choose the mesoscopic
scale $L'(t)=(\kappa t)^c$ and the fine scale $r=(\kappa t)^b$.
Note that by part $(ii)$ of lemma \ref{followinge},
 it is enough to prove that,

$$\frac{\sum_{x\in\Lambda_L}(\tilde m_a(x,t,w)-\langle\tilde m_a (x,t)\rangle)}
{\sqrt{Var_\mu\sum_{x\in\Lambda_L}\tilde m_a (x,t)}},
$$
converges in distribution to the normal law ${\mathcal N(0,1)}$. To do so,
we  write,

\begin{equation}
\label{ethat}
\frac{\sum_{x\in\Lambda_L}(\tilde m_a-\langle\tilde m_a\rangle)}
{\sqrt{Var_\mu\sum_{x\in\Lambda_L}\tilde m_a}}
=
\frac{\sum_{x\in S_L}(\tilde m_a-\langle\tilde m_a \rangle)}
{\sqrt{Var_\mu\sum_{x\in\Lambda_L}\tilde m_a }}
+\frac{\sum_{{\bf i}\in{\mathcal I }}
\sum_{x\in\Lambda''_{\bf i}}(\tilde m_a-\langle\tilde m_a\rangle)}
{\sqrt{Var_\mu\sum_{x\in\Lambda_L}\tilde m_a }}.
\end{equation}
We will first show that the strip component of the decomposition
(\ref{ethat}) converges to $0$ in probability. In fact, note
that for $t$ large enough, we have by statement (\ref{variance2})
of part $(ii)$ of lemma \ref{correlations} applied with
$U_t=S_{L(t)}$ and $U_t=\Lambda_{L(t)}$, that
$$
\frac{Var_\mu \sum_{x\in S_L}\tilde m_a}
{Var_\mu \sum_{x\in \Lambda_L}\tilde m_a}\sim
\frac{|S_L|}{(2L+1)^d}\le \min\{\kappa,(\kappa t)^{-(c-b)}\},
$$
where for the last inequality we have used estimate (\ref{p4}). Since $c>b$, this
tends to $0$ as $t\to\infty$.

Therefore, it is enough to prove that the second term of the
right-hand side of equality (\ref{ethat}), tends in law to
${\mathcal N}(0,1)$. For this purpose, since
the random variables $\left\{\sum_{x\in\Lambda''_{\bf i}}
(\tilde m_a-\langle\tilde m_a\rangle):{\bf i}\in{\mathcal I}\right\}$
are independent, it is enough to verify a version of
the Lyapunov condition. Namely, we will show that,

\begin{equation}
\label{owthat}
\frac{\sum_{{\bf i}\in{\mathcal I}}\left\langle\left|
\sum_{x\in\Lambda''_{\bf i}}(\tilde m_a -\langle\tilde m_a\rangle)
\right|^{2+\epsilon}\right\rangle}
{\left(\sum_{{\bf i}\in{\mathcal I}} Var_\mu\sum_{x\in\Lambda''_{\bf i}}
\tilde m_a\right)^{1+\epsilon/2}}\ll 1,
\end{equation}
and then apply again statement
(\ref{variance2})
of part $(ii)$ of lemma \ref{correlations}  to conclude that the variance
of the denominator
of the second term of display can be substituted by
$\sum_{{\bf i}\in{\mathcal I}} Var_\mu\sum_{x\in\Lambda''_{\bf i}}
\tilde m_a$.
Now, by the same token, we see that the denominator of the left-hand side of display (\ref{owthat}),
behaves asymptotically as $t\to\infty$ like

\begin{equation}
\label{like}
(2L+1)^{d(1+\epsilon/2)}\left(\sum_{x\in{\mathbb Z}^d}c_a(0,x,t)\right)^{1+\epsilon/2}.
\end{equation}
Thus, it is enough to prove the asymptotic behavior
 (\ref{owthat}), with the denominator
replaced by (\ref{like}). Now, note that
$\sum_{x\in{\mathbb Z}^d}c_a(0,x,t)$ is lower
bounded by $\langle \tilde m_a^2(0,t)\rangle-\langle\tilde m_a(0,t)\rangle^2$.
And by parts $(i)$ and $(iii)$ of proposition \ref{numbers} and lemma
\ref{improvement},
  this variance is lower-bounded by an expression which
is asymptotically equivalent to,

\begin{equation}
\label{vv}
c\left((\kappa t)^{3da}+1)^{-1}\langle m(0,2t)\rangle-
\langle m(0,t)\rangle^2\right).
\end{equation}
Now, assumption {\bf (MI)} with $\theta =2$ and
the inequality $\langle \tilde m_a(0,t)\rangle\le\langle m(0,t)\rangle$,
 imply that,
$(\log\langle m(0,2t)\rangle-2\log\langle m(0,t)\rangle\gg \log t$.
This shows that
 the second term $c\langle m(0,t)\rangle^2$ of the  expression
(\ref{vv}),
 is negligible with
respect to the first one $c ((\kappa t)^{3da}+1)^{-1}\langle m(0,2t)\rangle$.
 Hence,
 we conclude that it is enough to prove that,

\begin{equation}
\label{owthat2}
((\kappa t)^{3da}+1)
\frac{\sum_{{\bf i}\in{\mathcal I}}\left\langle\left|
\sum_{x\in\Lambda''_{\bf i}}(\tilde m_a -\langle\tilde m_a\rangle)
\right|^{2+\epsilon}\right\rangle}
{(2L+1)^{d(1+\epsilon/2)}
\langle m(0,2t)\rangle^{\left(1+\frac{\epsilon}{2}\right)}
}\ll 1.
\end{equation}
 By Jensen's inequality (\ref{jensen}) and
the upper-bound in display (\ref{p2}), we see that,
$\sum_{{\bf i}\in{\mathcal I}}\left\langle\left|
\sum_{x\in\Lambda''_{\bf i}}(\tilde m_a -\langle\tilde m_a\rangle)
\right|^{2+\epsilon}\right\rangle
\le (L'+1)^{d(1+\epsilon)}(2L+1)^d\left\langle\left|\tilde m_a-
\langle \tilde m_a\rangle\right|^{2+\epsilon}\right\rangle$.
Using again the fact that for
non-negative reals $x,y$ we have $|x-y|^{1+\epsilon}
\le |x|^{1+\epsilon}+|y|^{1+\epsilon}$, it follows that
$\left\langle\left|\tilde m_a-
\langle \tilde m_a\rangle\right|^{2+\epsilon}\right\rangle
\le \left\langle m^{2+\epsilon}\right\rangle+
\left\langle m\right\rangle^{2+\epsilon}
\le 2\left\langle m^{2+\epsilon}\right\rangle$, where we
have used Jensen's inequality in the last inequality. We see
that the left-hand side of (\ref{owthat}) is upper-bounded
by,

$$
((\kappa t)^{3da}+1)
(L'+1)^{d(1+\epsilon)}
\frac{\langle m(0,t)^{2+\epsilon}\rangle}{(2L+1)^{d\frac{\epsilon}{2}}
\langle m(0,2t)\rangle^{\left(1+\frac{\epsilon}{2}\right)}
}.
$$
Finally, by the hypothesis on the growth of $L$ this can be bounded
by,

$$((\kappa t)^{3da}+1)^2
(L'+1)^{d(1+\epsilon)}
e^{\frac{\epsilon}{2} F_{\frac{\epsilon}{2}}(2t)-
\frac{\epsilon}{2}F_\epsilon (2t)},
$$
which by condition {\bf (MI)},
 tends to $0$ as $t\to\infty$.

\medskip

\subsection{The non-Gaussian asymptotics}  Here we will prove
the asymptotics of display (\ref{int4}) of part $(ii)$ of theorem
\ref{te1}. It will be necessary to
perform a parity partition of $\Lambda_L$ with
a mesoscopic scale $L'=(\kappa t)^b$ for some $b>1$.
First note that  by display (\ref{item2}) of part $(ii)$ of lemma \ref{followinge} it will be enough to show that for some $b>a>1$,

\begin{equation}
\label{howthat}
\left\langle\left|
\frac{\sum_{x\in\Lambda_L}(\tilde m_a-\langle\tilde m_a\rangle)}
{\sqrt{Var_\mu\sum_{x\in\Lambda_L} m}}
\right|^{2-\epsilon}\right\rangle\ll 1.
\end{equation}
Now, by display (\ref{variance1}) of part $(ii)$
of lemma \ref{correlations},  the denominator
$\left(Var_\mu\sum_{x\in\Lambda_L} m\right)^{1-
\frac{\epsilon}{2}}$, of
the left-hand side of display (\ref{howthat}),
can be lower-bounded by
$(2L+1)^{d\left(1-\frac{\epsilon}{2}\right)}\left(\left\langle
 m^2\right\rangle-\langle m\rangle^2\right)^{1-\frac{\epsilon}{2}}\ge
c (2L+1)^{d\left(1-\frac{\epsilon}{2}\right)}\left\langle
 m^2\right\rangle^{1-\frac{\epsilon}{2}}
\ge
c ((\kappa t)^{3da}+1)^{-1}(2L+1)^{d\left(1-\frac{\epsilon}{2}\right)}\left\langle
 m(0,2t)\right\rangle^{1-\frac{\epsilon}{2}}$, where in the second to last
inequality
we used the assumption {\bf (MI)} and lemma \ref{improvement}
and in the last inequality assumption \ref{improvement}.
On the other hand, inequality  (\ref{truet}) of lemma \ref{ofand},
applied with $r=2-\epsilon$, Jensen's inequality and
$\tilde m_a\le m$, shows us that the numerator of
the left-hand side of display (\ref{howthat}) is upper-bounded
by $4(2L'+1)^{d(1-\epsilon)}(2L+1)^{d\frac{\epsilon}{2}}
 \langle m^{2-\epsilon}\rangle$.
Using the upper-bound $\langle m^{2-\epsilon}(0,t)\rangle\le
 k_1 ((\kappa t)^{3da}+1)\langle m(0,(2-\epsilon)t)\rangle$ of
lemma \ref{improvement}, the definition of the intermittency
       exponents $\{F_\theta\}$
in (\ref{ef}), and the assumption
$L(t)\le e^{\frac{1}{2}F_{-\epsilon}(2t)}$,
we hence see that it is enough to show that,

$$(2L'+1)^{d(1-\epsilon)}((\kappa t)^{3da}+1)^2
e^{\frac{\epsilon}{2}F_{-\epsilon}(2t)}
e^{-\frac{\epsilon}{2}F_{-\frac{\epsilon}{2}}(2t)}\ll 1.$$
But this is a consequence of assumption {\bf (MI)}.

\bigskip

\medskip

\subsection{Proof of theorem \ref{te2}}
Let us first prove part $(i)$
for Weibull-type tails. Note that it is enough to find an $x>0$
such that,

$$
\left\langle\left(
\frac{m_L(0,t)}{e^{(a(\gamma)+\delta)H(t)}}\right)^x\right\rangle\ll 1.
$$
Now, this expression is upper bounded by,

$$
e^{-(x(a+\delta)H(t)-H(xt)-(1-x)\gamma H(t))+o(H(t))},
$$
where we have used the asymptotics  (\ref{gm}). Since $H\in R_{\rho'}$
for $\rho'=\frac{\rho}{\rho-1}$, we see that the above expression is
upper bounded by,

$$
e^{-f_W(x,a+\delta) H(t)+o(H(t))},
$$
where $f_W(x,b):=x b-x^{\rho'}-(1-x)\gamma$ for $x>0, b>0$. This function has
a unique root at $x_0=\left(\frac{\rho-1}{\rho}(b+\gamma)\right)^{\rho-1}$
when $b=a(\gamma)$. Choosing $x=x_0$, we get the upper bound,

$$e^{-\delta x_0 H(t)+o(H(t))}\ll 1.$$
This proves part $(i)$ of theorem \ref{te2}. The proof of part $(iv)$
for Fr\'echet-type tails is completely analogous to the previous
argument, so it will be omitted. To prove part $(ii)$, note
that in analogy to the proof of part $(i)$, it is enough to show
that,

\begin{equation}
\label{tt1}
e^{-\left(x\frac{H((a+\delta)t)}{a+\delta}-H(xt)-(1-x)\gamma t\right)+o(t)}\ll 1.
\end{equation}
Now, by supposition  (\ref{ass-h}) we have
$\frac{H(xt)-xH((a+\delta)t)/(a+\delta)}{t}\sim -\rho x\log\frac{x}{a+\delta}.
$
Thus, the expression of display (\ref{tt1}), is upper-bounded by,

$$
e^{-f_D(x,a+\delta) t+o(t)},
$$
where $f_D(x,b):=\rho x\log\frac{x}{b}-(1-x)\gamma$ for $x>0, b>0$. But
this function has a single root at $x_0=ae^{\frac{1}{\rho}(\gamma-1)}$
when $b=a$. Hence, we obtain the upper bound,

$$
e^{-x_0\log\left(1+\frac{\delta}{a}\right)+o(t)}\ll 1.
$$
This proves part $(ii)$ of theorem \ref{te2}. The proof of part $(iii)$
for almost bounded potentials is analogous to the proof of part $(ii)$
so it will be omitted.

\medskip

\noindent  {\bf Acknowledgments.}
 Alejandro Ram\'\i rez
acknowledges
Fondecyt Grants 1020686 and 7020686 for their
financial support. He furthermore thanks the Courant
Institute of Mathematical Science, New York City, for
its kind hospitality, where part of this work was done.

\end{document}